\input amstex
\input epsf
\documentstyle{amsppt}\nologo\footline={}

\def\Gr{\mathop{\text{\rm Gr}}}
\def\tr{\mathop{\text{\rm tr}}}
\def\T{\mathop{\text{\rm T}}}
\def\Lin{\mathop{\text{\rm Lin}}}
\def\Re{\mathop{\text{\rm Re}}}
\def\E{\mathop{\text{\rm E}}}
\def\S{\mathop{\text{\rm S}}}
\def\B{\mathop{\text{\rm B}}}
\def\G{\mathop{\text{\rm G}}}
\def\D{\mathop{\text{\rm D}}}
\def\L{\mathop{\text{\rm L}}}
\def\ta{\mathop{\text{\rm ta}}}
\def\arccosh{\mathop{\text{\rm arccosh}}}
\def\Tn{\mathop{\text{\rm Tn}}}
\def\Ct{\mathop{\text{\rm Ct}}}
\def\Eu{\mathop{\text{\rm Eu}}}
\def\Area{\mathop{\text{\rm Area}}}

\hsize450pt

\topmatter\title Coordinate-free classic geometries\endtitle\author
Sasha Anan$'$in and Carlos H.~Grossi\endauthor\thanks First author
partially supported by the Institut des Hautes \'Etudes Scientifiques
(IH\'ES).\endthanks\thanks Second author partially supported by
FAPEMIG (grant 00163/06) and by the Max-Planck-Gesellschaft.\endthanks
\address Departamento de Matem\'atica, IMECC, Universidade Estadual de
Campinas,\newline13083-970--Campinas--SP, Brasil\endaddress\email
Ananin$_-$Sasha\@yahoo.com\endemail\address Max-Planck-Institut f\"ur
Mathematik, Vivatsgasse 7, 53111 Bonn, Germany\endaddress\email
grossi$_-$ferreira\@yahoo.com\endemail\subjclass 53A20 (53A35,
51M10)\endsubjclass\abstract This paper is devoted to a coordinate-free
approach to several classic geometries such as hyperbolic (real,
complex, quaternionic), elliptic (spherical, Fubini-Study), and
lorentzian (de Sitter, anti de Sitter) ones. These geometries carry a
certain simple structure that is in some sense stronger than the
riemannian structure. Their basic geometrical objects have linear
nature and provide natural compactifications of classic spaces. The
usual riemannian concepts are easily derivable from the strong
structure and thus gain their coordinate-free form. Many examples
illustrate fruitful features of the approach. The framework introduced
here has already been shown to be adequate for solving problems
concerning particular classic spaces.\endabstract\keywords Classic
geometries, hyperbolic geometry, Fubini-Study metric, parallel
transport\endkeywords\endtopmatter\document

\centerline{\bf1.~Classic geometries: introduction, definition,
examples, and motivation}

\bigskip

{\bf1.1.~Introduction.} The present paper constitutes an attempt to
systematically develop a coordinate-free view on several classic
geometries. The approach originates from [AGG] where, in order to
simplify formulae, we expressed several complex hyperbolic geometry
concepts in an invariant (hence, more convenient) form.

The riemannian structure in many classic geometries (hyperbolic,
spherical, Fubini-Study, etc.) turns out to be a shadow of a simpler
one. Let us briefly describe this stronger structure. Take a
$\Bbb K$-vector space $V$ with an hermitian form. The tangent vectors
to the grassmannian $\Gr_\Bbb K(r,V)$ at a nondegenerate point $p$ are
known to be $\Bbb K$-linear maps $p\to p^\perp$. We believe that a more
adequate object should be simply a $\Bbb K$-linear map $V\to V$, a {\it
footless\/} tangent vector: being composed with the two projectors
related to $p$, i.e., being {\it observed\/} from $p$, it becomes a
usual tangent vector. The {\it product\/} $t_1^*t_2$ (where
$t_1,t_2:V\to V$ are $\Bbb K$-linear maps and $t_1^*$ stands for the
map adjoint to $t_1$) is the structure that provides the hermitian
(riemannian) metric given by $\langle t_1,t_2\rangle:=\tr(t_1^*t_2)$
for $t_1,t_2$ observed from the same point~$p$. The
$(2,1)$-symmetrization of the triple product $tt_2^*t_1$ provides the
curvature tensor $R(t_2,t_1)t$ for $t,t_1,t_2$ observed from the same
point. Taking more observers in the previous examples, we obtain more
geometric characteristics. Distance, for instance, appears when one
observer sees the mote in the other observer's eye, i.e., when the
projectors related to their points are composed.

The basic objects in a classic geometry are linear in nature. This
makes grassmannians (and $\Bbb C$-grassmannians; see Subsection 1.7) a
place where these objects naturally vary. So, grassmannians should be
studied even if one is interested only in geometries embedded into
projective spaces. Regarding a classic geometry as a homogeneous space
related to the corresponding unitary (or orthogonal) group is
deficient: this does not allow to go outside the absolute, which would
be useful for the following reasons. The absolute (formed by degenerate
points) divides $\Gr_\Bbb K(r,V)$ into riemannian and pseudo-riemannian
pieces. Only one of them is traditionally considered as a classic
geometry. The grassmannian can be therefore seen as its
compactification. The points in each piece are in fact basic
geometrical objects (living in the traditional piece) whose type is
related to the compactification. Each piece is equipped with its
natural (pseudo-)riemannian geometry. Such geometries fit each other:
geometrical objects (geodesics, totally geodesic subspaces, equidistant
loci, etc.) pass through the absolute, leaving one piece
and entering another. Moreover, this global picture sheds light on the
geometry of the absolute. In particular, the general structure
described above (the one that provides the hermitian metric at
nondegenerate points) is inherited by the absolute. In the case of real
hyperbolic space, for instance, this explains the interrelation between
the conformal structure on the absolute and the metric structure on the
ball.

In classic geometries, the geometrical concepts and objects can be
introduced and handled synthetically. This suggests the above
modification of the usual riemannian tools and leads to simple linear
and hermitian algebra.

\smallskip

Some aspects of the coordinate-free approach can be found in
literature, including several examples of how such a framework was
successfully used in the solution of problems concerning particular
classic spaces. The following is a (very likely incomplete) list of
references:

\smallskip

$\bullet$ Concept of a projective model [Kle];

$\bullet$ Coordinate-free description of some particular metrics
[Arn1], [BeP];

$\bullet$ Linear approach to elementary geometric objects such as
geodesics, totally geodesic spaces, and bisectors [ChG], [Gir],
[Hsi1], [Hsi2], [Wil1], [Wil2];

$\bullet$ Linear and hermitian tools in real or complex hyperbolic
geometry [Gol], [HSa], [San], [Thu];

$\bullet$ Lorentzian projective compactification of real hyperbolic
space [Arn2], [ChK];

$\bullet$ Geometry of spaces of geodesics [AGK], [GeG], [GuK1], [Sal1],
[Sal2], [Stu];

$\bullet$ Solution of the Caratheodory conjecture [GuK2];

$\bullet$ Construction of new complex hyperbolic manifolds [AGG];

$\bullet$ Conformal structure on the absolute [AGoG];

$\bullet$ Strong structure on grassmannians [AGoG], [AGr].

\smallskip

In this article, we study projective classic geometries and describe in
a coordinate-free way several features of such geometries. In
particular, we obtain explicit expressions for the parallel
transport along geodesics in terms of the hermitian form
(Corollaries 5.7 and 5.9). Applying these expressions to the case of
complex hyperbolic geometry, we get a geometrical interpretation of the
angle between cotranchal bisectors in $\Bbb H^2_\Bbb C$ (Examples 6.1
and 6.3). Other explicit formulae involving geodesics (Subsections 3.2
and 3.4), projective cones (Example 3.6), bisectors (Examples 3.6 and
6.4), the Levi-Civita connection (Proposition 4.4), the curvature
tensor (Subsection 4.5), and sectional curvatures (Subsection 4.6) are
also provided.

For a similar treatment of grassmannian classic geometries, see [AGoG]
and [AGr].

\medskip

{\bf1.2.~Definition.} Let $\Bbb K$ denote one of the following fields:
$\Bbb R$ (real numbers), $\Bbb C$ (complex numbers), or~$\Bbb H$
(quaternions). A {\it classic geometry\/} is a right $\Bbb K$-vector
space $V$ equipped with an {\it hermitian form\/} $\langle-,-\rangle$.
By definition (see, for instance, [Lan]), the form is hermitian if it
takes values in $\Bbb K$, is biadditive, and satisfies the identities
$\langle v_1k,v_2\rangle=\overline k\langle v_1,v_2\rangle$,
$\langle v_1,v_2k\rangle=\langle v_1,v_2\rangle k$, and
$\langle v_1,v_2\rangle=\overline{\langle v_2,v_1\rangle}$ for all
$v_1,v_2\in V$ and $k\in\Bbb K$
$_\blacksquare$

\medskip

Behind this definition there is indeed more geometry than it might
appear at the first glance. The tangent space to a point $p$ in the
projective space $\Bbb P_\Bbb KV$ has a well-known description as the
$\Bbb R$-vector space ($\Bbb C$-vector space if $\Bbb K=\Bbb C$)
$${\T}_p\,\Bbb P_\Bbb KV={\Lin}_\Bbb K(p,V/p)\leqno{\bold{(1.3)}}$$
of all $\Bbb K$-linear transformations from $p$ to $V/p$. Here and in
what follows, we frequently do not distinguish the notation of a point
in $\Bbb P_\Bbb KV$, of a chosen representative of it in $V$, and of
the corresponding one-dimensional subspace when a concept or expression
does not depend on interpretation. For instance, the subspace $p^\perp$
of $V$ is well defined for any $p\in\Bbb P_\Bbb KV$.

If $p$ is nonisotropic, that is, if $\langle p,p\rangle\ne0$, then we
can naturally identify $V/p$ with $p^\perp$. In this case, we interpret
the tangent space as $\T_p\Bbb P_\Bbb KV=\Lin_\Bbb K(p,p^\perp)$. It
inherits the $\Bbb R$-bilinear form
$$(t_1,t_2):=\pm\frac{{\tr}_\Bbb R(t_1^*t_2)}{\dim_\Bbb R\Bbb
K},\leqno{\bold{(1.4)}}$$
where $t_1,t_2:p\to p^\perp$ are tangent vectors, $t_1^*:p^\perp\to p$
stands for the map adjoint to $t_1$ in the sense of the hermitian form,
and $\tr_\Bbb R(t_1^*t_2)$ denotes the trace of the $\Bbb R$-linear map
$t_1^*t_2:p\to p$. We will refer to this form as the {\it metric\/} of
a classic geometry. In the case of $\Bbb K=\Bbb C$, we have the {\it
hermitian metric}
$$\langle t_1,t_2\rangle:=\pm{\tr}_\Bbb
C(t_1^*t_2).\leqno{\bold{(1.5)}}$$
It is easy to see that $\Re\langle t_1,t_2\rangle=(t_1,t_2)$.
Obviously, the (hermitian) metric depends smoothly on a nonisotropic
$p$. If the hermitian form on $V$ is nondegenerate, then the metric is
nondegenerate. We~warn the reader that the case $\Bbb K=\Bbb H$
contains some peculiarities. The tangent space $\T_p\Bbb P_\Bbb HV$ is
not an $\Bbb H$-vector space and it makes no sense to speak of an
hermitian metric on it.

\smallskip

The {\it signature\/} of a point divides $\Bbb P_\Bbb KV$ into three
parts: {\it negative points,} {\it null points,} and {\it positive
points,} defined respectively as
$$\B V:=\{p\in\Bbb P_\Bbb KV\mid\langle p,p\rangle<0\},\qquad\S
V:=\{p\in\Bbb P_\Bbb KV\mid\langle p,p\rangle=0\},\qquad\E
V:=\{p\in\Bbb P_\Bbb KV\mid\langle p,p\rangle>0\}.$$

{\bf1.6.~Examples.} Take

\smallskip

(1)~$\Bbb K=\Bbb C$, $\dim_\Bbb CV=2$, the form of signature $++$,
and the sign $+$ in the definition of the hermitian metric. We obtain
the usual $2$-dimensional sphere of constant curvature.

\smallskip

(2)~$\Bbb K=\Bbb C$, $\dim_\Bbb CV=2$, the form of signature $+-$, and
the sign $-$ in the definition of the hermitian metric. Let
$p\in\Bbb P_\Bbb CV$ be nonisotropic. From the orthogonal decomposition
$V=p\oplus p^\perp$ it follows that the hermitian metric on
$\T_p\Bbb P_\Bbb CV$ is positive definite. We get two hyperbolic
Poincar\'e discs $\B V$ and $\E V$.

\smallskip

(3)~$\Bbb K=\Bbb R$, $\dim_\Bbb RV=3$, the form of signature $++-$, and
the sign $-$. We arrive at the hyperbolic Beltrami-Klein disc $\B V$.

\smallskip

(4)~$\Bbb K=\Bbb C$, $\dim_\Bbb CV=3$, the form of signature $++-$, and
the sign $-$. The open $4$-ball $\B V$ is the complex hyperbolic plane
$\Bbb H_\Bbb C^2$.

\smallskip

(5)~$\Bbb K=\Bbb H$, $\dim_\Bbb HV=2$, the form of signature $++$, and
the sign $+$. We obtain the usual $4$-sphere of constant curvature.
There is no $\Bbb H$-action on the tangent space $\T_p\Bbb P_\Bbb HV$.
However, fixing a geodesic in $\Bbb P_\Bbb HV$ leads to a curious
action of $\Bbb S^3\subset\Bbb H$ on the tangent bundle
$\T\Bbb P_\Bbb HV$ (see Example 3.7). The same is applicable to Example
1.6 (6) that follows.

\smallskip

(6)~$\Bbb K=\Bbb H$, $\dim_\Bbb HV=2$, the form of signature $+-$, and
the sign $-$. The open $4$-ball $\B V$ is the real hyperbolic space
$\Bbb H_\Bbb R^4$ (Example 3.7 shows a geometrical role of the
`additional' quaternionic structure).

\smallskip

In a similar way, we can describe many other geometries: elliptic
geometries such as spherical and Fubini-Study ones, hyperbolic
geometries including those of constant sectional or constant
holomorphic curvature, some lorentzian geometries such as de Sitter and
anti de Sitter spaces, etc.
$_\blacksquare$

\medskip

The most elementary geometrical objects are the `linear' ones, i.e.,
those given by the projectivization $\Bbb P_\Bbb KW$ of an
$\Bbb R$-vector subspace $W\subset V$. For instance, we can
isometrically embed Examples (1) and (2) as projective lines in Example
(4) by taking for $W$ an appropriate $2$-dimensional $\Bbb C$-vector
subspace in~$V$. (The negative part of a projective line of signature
$+-$ is commonly known as a {\it complex geodesic\/} in
$\Bbb H^2_\Bbb C$.) Let us take a look at some less immediate

\medskip

{\bf1.7. Examples.} (1)~Take $\dim_\Bbb RW=2$. Suppose that the
hermitian form, being restricted to $W$, is~real and does not vanish.
It is easy to see that $W\Bbb K\simeq W\otimes_\Bbb R\Bbb K$. The
circle
$$\G W:=\Bbb P_\Bbb KW=\Bbb P_\Bbb RW\simeq\Bbb S^1$$
is said to be a {\it geodesic.} The projective line
$\Bbb P_\Bbb K(W\Bbb K)$ is the {\it projective line of\/} the
geodesic. By Corollary~5.5, the introduced circle, out of its isotropic
points, is indeed a geodesic with respect to the metric and every
geodesic of the metric arises in this way.

\smallskip

(2)~Let $\dim_\Bbb RW=2$ in Example 1.6 (2). When $W$ is not a
$\Bbb C$-vector space (otherwise, $\Bbb P_\Bbb CW$ is simply a point in
$\Bbb P_\Bbb CV$), the real part of the hermitian form over $W$ can be
nondegenerate indefinite, definite, nonnull degenerate, or null. The
circle $\Bbb P_\Bbb CW$ is respectively said to be a {\it hypercycle,}
{\it metric circle,} {\it horocycle,} or the {\it absolute.} Inside
either of the Poincar\'e discs $\E V$ and $\B V$, we get the usual
hypercycles, metric circles, and horocycles.

\smallskip

(3)~We can isometrically embed (here the normalizing factor in (1.4)
plays its role) Example 1.6 (3) in Example 1.6 (4) by taking for $W$ a
$3$-dimensional $\Bbb R$-vector subspace such that the hermitian form,
being restricted to $W$, is real and nondegenerate. We obtain the {\it
$\Bbb R$-plane\/} $\Bbb P_\Bbb CW=\Bbb P_\Bbb RW\simeq\Bbb P_\Bbb R^2$,
a maximal lagrangian submanifold. The $\Bbb R$-planes are important in
complex hyperbolic geometry (see, for instance, [Gol] and [AGG]).

\smallskip

(4)~In Example 1.6 (4), let $S\subset V$ be an $\Bbb R$-vector
subspace, $\dim_\Bbb RS=2$. Suppose that the hermitian form is real and
nondegenerate over $S$. It is easy to see that $S^\perp$ is a
one-dimensional $\Bbb C$-vector space. Taking $W=S+S^\perp$, we arrive
at the {\it bisector\/} $B:=\Bbb P_\Bbb CW$. The geodesic $\G S$, the
projective line $\Bbb P_\Bbb C(S\Bbb C)$, and the point
$\Bbb P_\Bbb CS^\perp$ are respectively the {\it real spine,} the {\it
complex spine,} and the {\it focus\/} of the bisector. This description
of a bisector immediately provides (see [AGG, item 4.1.19]) the
well-known slice and meridional decompositions of a bisector (see
[Gir], [Mos], and [Gol]). If the hermitian form is indefinite over $S$,
then $\Bbb P_\Bbb CW\cap\B V$ is a usual bisector (= a hypersurface
equidistant from two points) in $\Bbb H^2_\Bbb C$. Every bisector in
$\Bbb H^2_\Bbb C$ is describable in this manner
$_\blacksquare$

\medskip

We would like to illustrate the thesis that the basic linear objects
form themselves spaces naturally endowed with a classic geometry
structure:

\smallskip

In Example 1.6 (3), the real projective plane $\Bbb P_\Bbb RV$ consists
of the usual Beltrami-Klein disc $\B V$ equipped with its riemannian
metric and of the M\"obius band $\E V$ endowed with a lorentzian
metric. The hermitian form establishes a duality between points and
projective lines ($=$ geodesics) in $\Bbb P_\Bbb RV$ : the point
$p\in\Bbb P_\Bbb RV$ corresponds to the geodesic
$\Bbb P_\Bbb Rp^\perp$. In view of this duality, the classic lorentzian
geometry of the M\"obius band $\E V$ is nothing but the geometry of
geodesics in the Beltrami-Klein disc $\B V$ and {\it vice versa.} By
the same reason, the classic pseudo-riemannian geometry of $\E V$ in
Example 1.6 (4) is the geometry of complex geodesics in
$\Bbb H^2_\Bbb C$.

\smallskip

In Example 1.7 (2), we will indistinctly refer to hypercycles, metric
circles, horocycles, and the absolute as {\it circles.} A point $W$ in
the grassmannian $\Gr_\Bbb R(2,V)$ of $2$-dimensional $\Bbb R$-vector
subspaces of $V$ determines in $\Bbb P_\Bbb CV$ a circle if $W$ is not
a $\Bbb C$-vector subspace and a point, otherwise. Clearly,
$W,W'\in\Gr_\Bbb R(2,V)$ provide the same circle if and only if $W=W'c$
for some $c\in\Bbb C^*$. The $\Bbb C$-{\it grassmannian\/}
$\Gr_{\Bbb C\mid\Bbb R}(2,V)$ is the quotient of $\Gr_\Bbb R(2,V)$ by
this action.

The singular locus of $\Gr_{\Bbb C\mid\Bbb R}(2,V)$ is formed by the
complex subspaces of $V$ and, therefore, coincides with
$\Bbb P_\Bbb CV$. It is easy to see that $\Gr_{\Bbb C\mid\Bbb R}(2,V)$
is topologically $\Bbb P_\Bbb R^3$ without an open $3$-ball. It has
$\Bbb P_\Bbb CV$ as its boundary. The absolute, a $2$-sphere with a
single double point, is formed by the horocycles and divides
$\Gr_{\Bbb C\mid\Bbb R}(2,V)$ into two parts.

How can we equip the $\Bbb C$-grassmannian
$\Gr_{\Bbb C\mid\Bbb R}(r,V)$ of $r$-dimensional $\Bbb R$-vector
subspaces of $V$ with a classic geometry structure? Let
$W\in\Gr_{\Bbb C\mid\Bbb R}(r,V)$ be a {\it nondegenerate\/} point,
that is, the real form $(-,-):=\Re\langle-,-\rangle$ is nondegenerate
over $W$. A tangent vector in $\T_W\Gr_{\Bbb C\mid\Bbb R}(r,V)$ is an
$\Bbb R$-linear map $t:W\to W^\perp$ such that $\tr_\Bbb R(\pi_Wit)=0$,
where the orthogonal $W^\perp$ is taken with respect to $(-,-)$ and
$\pi_W$ is the orthogonal projection onto $W$. The metric is given by
$(t_1,t_2):=\tr_\Bbb R(t_1^*t_2)$, where $t_1^*:W^\perp\to W$ is the
map adjoint to $t_1$ in the sense of $(-,-)$.

\bigskip

\newpage

\centerline{\bf2.~Preliminaries}

\medskip

Let $p\in\Bbb P_\Bbb KV$ be a nonisotropic point. We introduce the
following notation of orthogonal decomposition
$$V=p\oplus p^\perp,\qquad v=\pi'[p]v+\pi[p]v,$$
where
$$\pi'[p]v:=p\,\frac{\langle p,v\rangle}{\langle p,p\rangle}\in p\Bbb
K\qquad\text{and}\qquad\pi[p]v:=v-p\,\frac{\langle p,v\rangle}{\langle
p,p\rangle}\in p^\perp$$
do not depend on the choice of a representative of $p$. Depending on
circumstances, we choose the most convenient variant of notation.

\smallskip

The hermitian form over a 2-dimensional $\Bbb K$-vector subspace of $V$
can be null, definite, nondegenerate indefinite, or nonnull degenerate.
The corresponding projective line will be respectively called {\it
null,} {\it spherical,} {\it hyperbolic,} or {\it euclidean.} We need a
very rudimental form of Sylvester's criterion applicable to the case
$\Bbb K=\Bbb H$.

\medskip

{\bf2.1.~Lemma.} {\sl Let\/ $W$ be a\/ $2$-dimensional\/
$\Bbb K$-vector space equipped with a nonnull hermitian form.
The hermitian form is respectively definite, nondegenerate indefinite,
or degenerate if and only if\/ $\D(p,q)>0$, $\D(p,q)<0$, or\/
$\D(p,q)=0$, where\/
$\D(p,q):=\langle p,p\rangle\langle q,q\rangle-\langle
p,q\rangle\langle q,p\rangle$
and\/ $p,q$ are any two\/ $\Bbb K$-linearly independent vectors in\/
$W$. {\rm(}Obviously, $\D(p,q)=0$ if\/ $p,q$ are\/ $\Bbb K$-linearly
dependent.{\rm)}}

\medskip

{\bf Proof.} If one of $p,q$ is nonisotropic (say, $p$) the result
follows from $\pi[p]q\ne0$, $\big\langle p,\pi[p]q\big\rangle=0$, and
$$\langle p,p\rangle\big\langle\pi[p]q,\pi[p]q\big\rangle=\langle
p,p\rangle\big\langle q,\pi[p]q\big\rangle=\langle
p,p\rangle\Big(\langle q,q\rangle-\frac{\langle q,p\rangle\langle
p,q\rangle}{\langle p,p\rangle}\Big)=\D(p,q).$$
If both $p,q$ are isotropic, we take a nonisotropic $u\in W$. We can
assume that $u=pk+q$ for some $k\in\Bbb K^*$. Clearly, $\pi[u]q\ne0$,
$\big\langle u,\pi[u]q\big\rangle=0$, and
$\D(u,q)=\langle u,u\rangle\big\langle\pi[u]q,\pi[u]q\big\rangle$. It
remains to observe that
$$\D(u,q)=\langle pk+q,pk+q\rangle\langle q,q\rangle-\langle
pk+q,q\rangle\langle q,pk+q\rangle=-\overline k\langle
p,q\rangle\langle q,p\rangle k=|k|^2\D(p,q)\ _\blacksquare$$

{\bf2.2.~Remark.} (1)~Let $\L$ be a projective line. For every
nonisotropic $p\in\L$ there exists a unique $q\in\L$ {\it orthogonal\/}
to $p$, that is, such that $\langle p,q\rangle=0$.

\smallskip

(2)~Isotropic points in a hyperbolic projective line form an
$(n-1)$-sphere, where $n=\dim_\Bbb R\Bbb K$. An~euclidean projective
line contains a single isotropic point
$_\blacksquare$

\medskip

A linear transformation in (1.3) can be regarded as a tangent vector in
usual differential terms: Let $f$ be a $\Bbb K$-valued smooth function
defined in a neighbourhood of $p\in\Bbb P_\Bbb KV$ and let $\hat f$
denote its lift to the corresponding neighbourhood of
$p\Bbb K\setminus\{0\}$ in $V$. Clearly, $\hat f(vk)=\hat f(v)$ for all
$k\in\Bbb K^*$. Every $\varphi\in\Lin_\Bbb K(p,V)$ defines a tangent
vector $t_\varphi\in\T_p\Bbb P_\Bbb KV$ given by
$$t_{\varphi}f:=\frac d{d\varepsilon}\Big|_{\varepsilon=0}\hat
f\big((1+\varepsilon\varphi)p\big),$$
where $\varepsilon\in\Bbb R$. Note that $t_\varphi$ vanishes if and
only if $\varphi p\in p\Bbb K$. Also, altering $\varphi$ by adding $pk$
to $\varphi p$, where $k\in\Bbb K$, does not change the vector
$t_\varphi\in\Lin(p,V/p)$.

If $p\in\Bbb P_\Bbb KV$ is nonisotropic, we have the identification
$${\T}_p\Bbb P_\Bbb KV={\Lin}_\Bbb K(p,p^\perp)=p^\perp\langle
p,-\rangle.\leqno{\bold{(2.3)}}$$

\smallskip

{\bf2.4.~Remark.} (1)~In terms of (2.3), the map adjoint to
$v\langle p,-\rangle$ is given by
$\big(v\langle p,-\rangle\big)^*=p\langle v,-\rangle$, where
$v\in p^\perp$.

\smallskip

(2)~Let $p\in\Bbb P_\Bbb KV$ be nonisotropic and let $v\in V$. Then the
trace of the $\Bbb R$-linear map $t:=v\langle p,-\rangle$ is given by
$\tr_\Bbb Rt=\dim_\Bbb R\Bbb K\cdot\Re\langle p,v\rangle$.

\smallskip

This treatment is useful while performing explicit calculations
$_\blacksquare$

\medskip

{\bf2.5.~Definition.} Let $W\subset V$ be an $\Bbb R$-vector subspace.
We call a point $p\in W$ {\it projectively smooth\/} in $W$ if
$\dim_\Bbb R(p\Bbb K\cap W)=\min\limits_{0\ne w\in W}\dim_\Bbb R(w\Bbb
K\cap W)$
$_\blacksquare$

\medskip

It is not difficult to see that the projectively smooth points in $W$
provide an open smooth region in $\Bbb P_\Bbb KW$. Moreover, we have
the following

\medskip

{\bf2.6.~Lemma {\rm[AGG, Lemma 4.2.2]}.} {\sl Let\/ $W\subset V$ be
an\/ $\Bbb R$-vector subspace, let\/ $p\in W$ be a projectively smooth
point in\/ $W$, and let\/ $\varphi\in\Lin_\Bbb K(p,V)$. Then
$t_\varphi\in\T_p\Bbb P_\Bbb KW$ if and only if\/
$\varphi p\in W+p\Bbb K$}
$_\blacksquare$

\medskip

The tangent vector to a smooth path can be expressed in terms of the
identification $\T_p\Bbb P_\Bbb KV=\Lin_\Bbb K(p,p^\perp)$ :

\medskip

{\bf2.7.~Lemma {\rm[AGG, Lemma 4.1.4]}.} {\sl Let\/
$c:[a,b]\to\Bbb P_\Bbb KV$ be a smooth curve and let\/ $c_0:[a,b]\to V$
be a smooth lift of\/ $c$ to\/ $V$. If\/ $c(t_0)$ is nonisotropic, then
the tangent vector\/ $\dot c(t_0):c_0(t_0)\to c_0(t_0)^\perp$ is given
by\/ $\dot c(t_0):c_0(t_0)\mapsto\pi\big[c(t_0)\big]\dot c_0(t_0)$}
$_\blacksquare$

\bigskip

\centerline{\bf3.~Geodesics}

\medskip

Let us remind the definition in Example 1.7 (1). Take a $2$-dimensional
$\Bbb R$-vector subspace $W\subset V$ such that the hermitian form,
being restricted to $W$, is real and does not vanish. It is immediate
that $W\Bbb K\simeq W\otimes_\Bbb R\Bbb K$. Hence,
$\Bbb P_\Bbb KW=\Bbb P_\Bbb RW$. The circle $\G W:=\Bbb P_\Bbb KW$ is,
by definition, a geodesic. (Corollary~5.5 relates this concept to the
common one.) The geodesic $\G W$ spans its projective line
$\Bbb P_\Bbb K(W\Bbb K)$. A geodesic is called {\it spherical,} {\it
hyperbolic,} or {\it euclidean\/} depending on the nature of its
projective line.

\medskip

{\bf3.1.~Lemma.} {\sl{\rm(1)}~Let\/ $g_1,g_2\in\Bbb P_\Bbb KV$ be
distinct and nonorthogonal. Then there exists a unique geodesic
containing\/ $g_1$ and\/ $g_2$.

{\rm(2)}~Let\/ $p\in\Bbb P_\Bbb KV$ be nonisotropic and let\/
$0\ne t\in\T_p\Bbb P_\Bbb KV$, $t:p\to p^\perp$. Then there exists a
unique geodesic having\/ $t$ as its tangent vector at\/ $p$. This
geodesic is given by the subspace\/ $W=p\Bbb R+tp\Bbb R$.}

\medskip

{\bf Proof.} (1)~Clearly, $g_1,g_2\in\G W$ for
$W=g_1\Bbb R+g_2\langle g_2,g_1\rangle\Bbb R$. If $g_1,g_2\in\G W'$,
then $W'=g_2k_2\Bbb R+g_1k_1\Bbb R$ for some $k_1,k_2\in\Bbb K$ such
that $\overline k_2\langle g_2,g_1\rangle k_1\in\Bbb R^*$. Hence,
$W'=g_2k_2\overline k_2\langle g_2,g_1\rangle k_1\Bbb R+g_1k_1\Bbb
R=Wk_1$,
that is, $\G W'=\G W$.

(2)~The geodesic $\G W$, where $W=p\Bbb R+tp\Bbb R$, does not depend on
the choice of $p\in p\Bbb K$. By Lemma~2.6,  $t$ is a tangent vector to
$\G W$ at $p$. Let $\G W'$ be a geodesic with tangent vector $t$. We
can choose $W'$ so that $p\in W'$. By Lemma 2.6, $tp\in W'+p\Bbb K$.
So, $tp\in p^\perp$ implies $tp\in W'$. In other words,
$W'=p\Bbb R+tp\Bbb R$~
$_\blacksquare$

\medskip

We denote by $\G{\wr}g_1,g_2{\wr}$ the geodesic that contains given
distinct nonorthogonal $g_1,g_2\in\Bbb P_\Bbb KV$.

\smallskip

Take distinct {\bf orthogonal} $g_1,g_2\in\Bbb P_\Bbb KV$. Assume that
the projective line $\L$ spanned by $g_1,g_2$ is nonnull. One of
$g_1,g_2$ is nonisotropic --- say, $g_1$. Every geodesic in $\L$
passing through $g_1$ has the form $\G W$ with $W=q\Bbb R+g_1\Bbb R$,
$g_1\ne q\in\L$, and $\langle q,g_1\rangle\in\Bbb R^*$. So,
$\pi[g_1]q\in\G W$. By Remark 2.2 (1), $g_2$ is the only point in $\L$
orthogonal to $g_1$. Hence, $\pi[g_1]q=g_2$ in $\Bbb P_\Bbb KV$. In
other words, every geodesic in $\L$ that passes through $g_1$ also
passes through $g_2$. In particular, every geodesic in an euclidean
projective line passes through the isotropic point (see Remark 2.2
(2)). In this case, in the affine chart $\Bbb K$ of nonisotropic points
of $\L$, the geodesics correspond to the straight lines. This justifies
the term `euclidean.' Since the metric is actually null over euclidean
lines, perhaps a more appropriate term would be {\it affine line.}

\medskip

{\bf3.2.~Length of noneuclidean geodesics.} Take a spherical projective
line $\L$, take a point $g_1\in\L$, and choose the sign $+$ in the
definition (1.4) of the metric. Let $g'_1\in\L$ denote the point
orthogonal to $g_1$. Fixing representatives $g_1,g'_1\in V$ such that
$\langle g_1,g_1\rangle=\langle g'_1,g'_1\rangle=1$, we parameterize a
lift $c_0(t):=g_1\cos t+g'_1\sin t$ to $V$ of a segment of geodesic
$c=c(t)$ joining $g_1$ and $g_2:=c(a)$, where $t\in[0,a]$ and
$a\in[0,\pi/2]$. Since $\big\langle\dot c_0(t),c_0(t)\big\rangle=0$ and
$\big\langle c_0(t),c_0(t)\big\rangle=1$, it follows from Lemma 2.7
that $\big(\dot c(t),\dot c(t)\big)=1.$ Hence,
$\ell\,c=\displaystyle\int_0^a\sqrt{\big(\dot c(t),\dot c(t)\big)}=a$.
Noting that $\ta(g_1,g_2)=\cos^2a$, where
$$\ta(g_1,g_2):=\frac{\langle g_1,g_2\rangle\langle
g_2,g_1\rangle}{\langle g_1,g_1\rangle\langle
g_2,g_2\rangle},\leqno{\bold{(3.3)}}$$
we obtain
$$\ell\,c=\arccos\sqrt{\ta(g_1,g_2)}.$$
It follows immediately from Lemma 2.1 that, being $\L$ spherical,
$0\le\ta(g_1,g_2)\le1$. The first equality occurs exactly when
$g_1,g_2$ are orthogonal and the second, exactly when $g_1=g_2$.

If $\L$ is a hyperbolic projective line, similar arguments involving
$\cosh$, $\sinh$, and the sign $-$ for the metric show that the length
of a segment of geodesic $c$ that contains no isotropic points and
joins $g_1,g_2\in\L$ is given by
$$\ell\,c=\arccosh\sqrt{\ta(g_1,g_2)}.$$

In both cases, the distance is a monotonic function of the {\it
tance\/} $\ta(g_1,g_2)$ (see also [AGG,
Corollary~4.1.18])
$_\blacksquare$

\medskip

{\bf3.4.~Equations of a geodesic.} Let the geodesic
$\G{\wr}g_1,g_2{\wr}$ be noneuclidean and let $\L$ denote its
projective line. We will show that $x\in\L$ belongs to
$\G{\wr}g_1,g_2{\wr}$ if and only if
$$b(x,g_1,g_2):=\langle x,g_1\rangle\langle g_1,g_2\rangle\langle
g_2,x\rangle-\langle x,g_2\rangle\langle g_2,g_1\rangle\langle
g_1,x\rangle=0.$$
The proof is straightforward. The above equation does not depend on the
choice of representatives $x,g_1,g_2\in V$. If
$x\in\G{\wr}g_1,g_2{\wr}$, then $b(x,g_1,g_2)=0$ since the hermitian
form is real over $W$ and we can assume $x,g_1,g_2\in W$. Suppose that
$b(x,g_1,g_2)=0$ for some $x\in\L$. We can take $g_1,g_2\in W$ and
$x=g_1k+g_2$ for some $k\in\Bbb K$. The condition $b(x,g_1,g_2)=0$ is
equivalent to
$\big(\langle g_1,g_2\rangle\langle g_2,g_1\rangle-\langle
g_1,g_1\rangle\langle g_2,g_2\rangle\big)(k-\overline k)=0$.
Since $\L$ is not euclidean, we conclude from Lemma 2.1 that
$k\in\Bbb R$, that is, $x\in W$.

\smallskip

Let $g\in\G{\wr}g_1,g_2{\wr}$ and let
$\varphi\in\Lin_\Bbb K(g,V)$ be such that $t_\varphi\in\T_g\L$.
We will show that $t_\varphi\in\T_g\G{\wr}g_1,g_2{\wr}$ if and only if
$$t(\varphi g,g,g_1,g_2):=\langle\varphi g,g_1\rangle\langle
g_1,g_2\rangle\langle g_2,g\rangle+\langle g,g_1\rangle\langle
g_1,g_2\rangle\langle g_2,\varphi g\big\rangle-$$
$$-\langle\varphi g,g_2\rangle\langle g_2,g_1\rangle\langle
g_1,g\rangle-\langle g,g_2\rangle\langle g_2,g_1\rangle\langle
g_1,\varphi g\rangle=0.$$
It follows from $b(g,g_1,g_2)=0$ that
$$t(\varphi g+gk,g,g_1,g_2\big)=t(\varphi g,g,g_1,g_2)+\overline k\cdot
b(g,g_1,g_2)+b(g,g_1,g_2)\cdot k=t(\varphi
g,g,g_1,g_2)\leqno{\bold{(3.5)}}$$
for every $k\in\Bbb K$. Also, the equation $t(\varphi g,g,g_1,g_2)=0$
does not depend on the choice of representatives for $g,g_1,g_2$. We
take $g,g_1,g_2\in W$. If $t_\varphi\in\T_g\G{\wr}g_1,g_2{\wr}$, then
$\varphi g\in W+g\Bbb K$ by Lemma 2.6. Due to (3.5), we can assume
that $\varphi g\in W$. Hence, $t(\varphi g,g,g_1,g_2)=0$. Conversely,
suppose that $t(\varphi g,g,g_1,g_2\big)=0$. We can take $g=g_1r+g_2$
for some $r\in\Bbb R$ (interchanging $g_1$ and $g_2$ if necessary).
Since $t_\varphi\in\T_g\L$, it~follows from Lemma 2.6 that
$\varphi(g)=g_1k_1+g_2k_2$ for some $k_1,k_2\in\Bbb K$. Due to (3.5),
we can assume that $\varphi g=g_1k$. Now, the condition
$t(\varphi g,g,g_1,g_2)=0$ means that
$\big(\langle g_1,g_2\rangle\langle g_2,g_1\rangle-\langle
g_1,g_1\rangle\langle g_2,g_2\rangle\big)(k-\overline k)=0$.
By Lemma 2.1, $k\in\Bbb R$, that is, $\varphi g\in W$
$_\blacksquare$

\medskip

{\bf3.6.~Example: equations of the cone over a geodesic.} We take
$\dim_\Bbb KV=3$ and a nondegenerate hermitian form
$\langle-,-\rangle$. The hermitian form establishes a correspondence
between points and projective lines in $\Bbb P_\Bbb KV$ : the point
$p\in\Bbb P_\Bbb KV$ corresponds to the projective line
$\Bbb P_\Bbb Kp^\perp$. We call $p$ the {\it polar point\/} to
$\Bbb P_\Bbb Kp^\perp$.

\smallskip

Let $\G S$ be a noneuclidean geodesic. Clearly, $S^\perp$ is a
$\Bbb K$-vector space and $p:=\Bbb P_\Bbb KS^\perp$ is the
(nonisotropic, by Lemma 2.1) polar point to the projective line of
$\G S$. Therefore, $C:=\Bbb P_\Bbb K(S+S^\perp)$ is the projective cone
over $\G S$ with vertex $p$. All elements in $S+S^\perp$, except those
in $S^\perp$, are projectively smooth (see Definition 2.5).

A point $x\in\Bbb P_\Bbb KV$ that is different from $p$ belongs to $C$
if and only if $\pi[p]x\in\G S$. Hence, $x\in C$ means that
$b\big(\pi[p]x,g_1,g_2\big)=0$ (see Subsection 3.4), where
$g_1,g_2\in\G S$ are distinct nonorthogonal points. This implies that
$C$ is given by the equation
$$b(x,g_1,g_2)=0.$$

Let $c\in C$ be different from $p$ and let
$\varphi\in\Lin_\Bbb K(c,V)$. Define a linear map
$\psi\in\Lin_\Bbb K(g,V)$ by putting $g:=\pi[p]c$ and
$\psi g:=\pi[p]\varphi c$. Fix a representative $c\in S+S^\perp$.
Clearly, $g\in S$. If $t_{\varphi}\in\T_cC$, then
$\varphi c\in S+S^\perp+c\Bbb K$ by Lemma 2.6. This implies that
$\psi g\in S+g\Bbb K$, that is, $t_{\psi}\in\T_g\G S$. Conversely,
if~$t_{\psi}\in\T_g\G S$, then
$\psi g\in S+g\Bbb K\subset S+S^\perp+c\Bbb K$. Hence,
$\varphi c\in S+S^\perp+c\Bbb K$. In other words,
$t_\varphi\in\T_c\Bbb P_\Bbb KV$ is tangent to $C$ if and only if
$t\big(\pi[p]\varphi c,g,g_1,g_2\big)=0$, where $g_1,g_2$ are distinct
nonorthogonal points in $\G S$. This is equivalent to
$$t(\varphi c,c,g_1,g_2)=0.$$

In the case of $\Bbb K=\Bbb C$, the projective cone $C$ is nothing but
the bisector with the real spine $\G S$ (see Example 1.7 (4) and the
references therein). From the equation of the tangent space to a point
in a bisector, one derives the expression
$$n(q,g_1,g_2)=\Big(g_1\frac{\langle g_2,q\rangle}{\langle
g_2,g_1\rangle}-g_2\frac{\langle g_1,q\rangle}{\langle
g_1,g_2\rangle}\Big)i\langle q,-\rangle$$
of the normal vector $n(q,g_1,g_2)$ at $q$ to the bisector whose real
spine is $\G{\wr}g_1,g_2{\wr}$ (see (2.3) and [AGG, Proposition
4.2.11]). This last expression permits to calculate, in terms of the
hermitian form, the oriented angle between two bisectors with a common
slice (see [AGG, Lemma 4.3.1] and Example~6.3)~
$_\blacksquare$

\medskip

{\bf3.7.~Example: actions on tangent bundle given by the choice of a
geodesic.} We consider the case $\Bbb K=\Bbb H$. The tangent space to a
point in $\Bbb P_\Bbb HV$ is not an $\Bbb H$-vector space. In order to
define an action of the sphere $\Bbb S^3\subset\Bbb H$ over the tangent
bundle $\T\Bbb P_\Bbb HV$, we assume that $V$ is an
$(\Bbb H,\Bbb H)$-bimodule.

Let $p\in\Bbb P_\Bbb HV$ and let $\varphi\in\Lin_\Bbb H(p,V)$. Given
$k\in\Bbb S^3\subset\Bbb H$, we define the linear map
$k\varphi\in\Lin_\Bbb H(kp,V)$ by putting
$(k\varphi)(kx):=k(\varphi x)$ for all $x\in p$. In this way, we arrive
at the left action $(p,t_\varphi)\mapsto(kp,t_{k\varphi})$ of
$\Bbb S^3$ over the tangent bundle $\T\Bbb P_\Bbb HV$ (note that
changing $\varphi p$ by $\varphi p+pk'$ results in the same
$t_{k \varphi}$). It is easy to verify that $t_{k\varphi}$ is also the
image of $t_\varphi$ under the differential $d(k\cdot)_p$, where
$k\cdot:\Bbb P_\Bbb HV\to\Bbb P_\Bbb HV$ is induced by
$k\cdot:v\mapsto kv$.

Suppose that the $(\Bbb H,\Bbb H)$-bimodule structure is {\it
compatible\/} with the hermitian form, that is,
$\langle v_1,kv_2\rangle=\langle\overline kv_1,v_2\rangle$ for all
$v_1,v_2\in V$ and $k\in\Bbb H$. Then, for a nonisotropic $p$ and for
$t:p\to p^\perp$, we have $d(k\cdot)_p\,t=kt:kp\to(kp)^\perp$. Hence,
out of isotropic points, $k\cdot$ is an isometry.

It is well known that every $(\Bbb H,\Bbb H)$-bimodule has the form
$V=W\otimes_\Bbb R\Bbb H$, where
$W:=\{v\in V\mid kv=vk\text{ for every }k\in\Bbb H\}$ is the centre of
the bimodule. The bimodule structure is compatible with
$\langle-,-\rangle$ if and only if the form restricted to $W$ is real.
In other words, the choice of a bimodule structure compatible with the
hermitian form is equivalent to the choice of a linear geometrical
object $\Bbb P_\Bbb KW$ corresponding to a maximal real subspace $W$ in
$V$.

In the particular case of $\dim_\Bbb HV=2$, we get an action of
$\Bbb S^3$ over $\Bbb P_\Bbb HV$ by isometries that is determined by
the choice of an arbitrary geodesic $\G$. This geodesic is the
fixed-point set of the action. The orbit of every other point is a
$2$-sphere. Thus, we obtain some foliation of
$\Bbb P_\Bbb HV\setminus\G$ by $2$-spheres. The actions over
$\Bbb P_\Bbb HV$ for hyperbolic (Example 1.6 (6)) and elliptic (Example
1.6 (5)) geometries produce topologically distinct foliations
$_\blacksquare$

\medskip

In Section 5, we show that the geodesics introduced in Example 1.7 (1)
are indeed geodesics with respect to the metric, out of their isotropic
points. Thus, for the classic geometries, we can forget about the
variational characterization of geodesics and deal only with the
`linear' one, which is much easier.

\bigskip

\centerline{\bf4.~Levi-Civita connection}

\medskip

From now on, we assume the hermitian form $\langle-,-\rangle$ to be
nondegenerate. In particular, $\B V$ and $\E V$ are endowed with
pseudo-riemannian metrics.

Also, until the end of the article, we use the following conventions.
Let $p\in\Bbb P_\Bbb KV$ be nonisotropic. Extending by zero, we
consider any tangent vector $t:p\to p^\perp$ as a linear map
$t\in\Lin_\Bbb K(V,V)$.
So,~$\T_p\Bbb P_\Bbb KV=\Lin_\Bbb K(p,p^\perp)\subset\Lin_\Bbb K(V,V)$.
(Obviously, $t=t\pi'[p]$, $t=\pi[p]t$, $t\pi[p]=\pi'[p]t=0$, and $st=0$
for all tangent vectors $s,t\in\Lin_\Bbb K(V,V)$ at $p$.) Conversely,
given an arbitrary linear map $t\in\Lin_\Bbb K(V,V)$, we define the
tangent vector
$$t_p:=\pi[p]t\pi'[p]$$
at $p$.

\smallskip

Let $U\subset V$ be a {\it saturated\/} open set (i.e.,
$U\Bbb K^*\subset U$) without isotropic points. A {\it lifted field\/}
over $U$ is a smooth map $X:U\to\Lin_\Bbb K(V,V)$ such that
$X(p)_p=X(p)$ and $X(pk)=X(p)$ for all $p\in U$ and $k\in\Bbb K^*$. In
other words, $X$ correctly defines a smooth tangent field over
$\Bbb P_\Bbb KU$.

\medskip

{\bf4.1.~Definition.} Every $t\in\Lin_\Bbb K(V,V)$ provides the
(lifted) field $T$ {\it spread\/} from $t$ : it is given by the rule
$T(p)=t_p$ and is defined for all nonisotropic $p$
$_\blacksquare$

\medskip

For $t\in\Lin_\Bbb K(V,V)$, we put
$$\nabla_tX(p):=\Big(\displaystyle\frac
d{d\varepsilon}\Big|_{\varepsilon=0}X\big((1+\varepsilon
t)p\big)\Big)_p.$$
Since $\pi[pk]=\pi[p]$ and $\pi'[pk]=\pi'[p]$ for all $p\in U$ and
$k\in\Bbb K^*$, the field $p\mapsto\nabla_{Y(p)}X$ is lifted for
arbitrary lifted fields $X$ and $Y$ over $U$. Obviously, $\nabla$
enjoys the properties of an affine connection.

\medskip

{\bf4.2.~Lemma.} {\sl Let\/ $p\in\Bbb P_\Bbb KV$ be nonisotropic and
let\/ $t$ be a tangent vector at\/ $p$. Then
$$\frac{d}{d\varepsilon}\Big|_{\varepsilon=0}\pi'[p+tp\varepsilon]=
-\frac{d}{d\varepsilon}\Big|_{\varepsilon=0}\pi[p+tp\varepsilon]=
t+t^*.$$}

{\bf Proof.} By definition,
$\displaystyle\pi'[p+tp\varepsilon]=(p+tp\varepsilon)\frac{\langle
p+tp\varepsilon,-\rangle}{\langle p,p\rangle+\varepsilon^2\langle
tp,tp\rangle}$.
Differentiating, we get
$$\frac{d}{d\varepsilon}\Big|_{\varepsilon=0}(p+tp\varepsilon)
\frac{\langle p+tp\varepsilon,-\rangle}{\langle
p,p\rangle+\varepsilon^2\langle tp,tp\rangle}=p\frac{\langle
tp,-\rangle}{\langle p,p\rangle}+tp\frac{\langle p,-\rangle}{\langle
p,p\rangle}.$$
The second term equals $t\pi'[p]=t$. Put
$\varphi:=p\displaystyle\frac{\langle tp,-\rangle}{\langle
p,p\rangle}$. Then
$$\langle tx,y\rangle=\big\langle t\pi'[p]x,y\big\rangle=\Big\langle
tp\frac{\langle p,x\rangle}{\langle
p,p\rangle},y\Big\rangle=\frac{\langle x,p\rangle}{\langle
p,p\rangle}\langle tp,y\rangle=\Big\langle x,p\frac{\langle
tp,y\rangle}{\langle p,p\rangle}\Big\rangle=\langle x,\varphi
y\rangle$$
for every $x,y\in V$. Hence, $t^*=\varphi$
$_\blacksquare$

\medskip

{\bf4.3.~Lemma.} {\sl Let\/ $p\in\Bbb P_\Bbb KV$ be nonisotropic. Let\/
$s$ and\/ $t$ be tangent vectors at\/ $p$. Then
$$\nabla_{T}S(x)=\big(s\pi[x]t-t\pi'[x]s\big)_x$$
for every nonisotropic\/ $x\in\Bbb P_\Bbb KV$, where the fields\/ $S$
and\/ $T$ are respectively\/ spread from\/ $s$ and\/ $t$\/ {\rm(}see
Definition\/ {\rm4.1).} In particular, $\nabla_{T}S(p)=0$.}

\medskip

{\bf Proof.} By Lemma 4.2,
$$\nabla_{T}S(x)=\nabla_{t_x}S(x)=\Big(\frac
d{d\varepsilon}\Big|_{\varepsilon=0}S(x+t_xx\varepsilon)\Big)_x=
\Big(\frac d{d\varepsilon}\Big|_{\varepsilon=0}\pi[x+t_xx\varepsilon]s
\pi'[x+t_xx\varepsilon]\Big)_x=$$
$$=\Big(-\big(t_x+(t_x)^*\big)s\pi'[x]+
\pi[x]s\big(t_x+(t_x)^*\big)\Big)_x=\big(s\pi[x]t-t\pi'[x]s\big)_x$$
since $\pi[x](t_x)^*=(t_x)^*\pi'[x]=0$
$_\blacksquare$

\medskip

The fact that $\nabla$ is Levi-Civita for the (hermitian) metric can be
easiy inferred from the theory of classical groups. Indeed, one needs
essentially to show that $\nabla$ is torsion-free and this holds
because there are no $3$-tensors which are invariant under the
orthogonal, unitary, or symplectic groups; see [Wey] (or~[How] for a
more modern treatment). However, we found it helpful to present below a
straightforward proof of the fact in question as it may illustrate the
role of spread fields (see Definition 4.1) and keep the exposition more
self-contained.

\medskip

{\bf4.4.~Proposition.} {\sl$\nabla$ is the Levi-Civita connection for
the\/ {\rm(}hermitian\/{\rm)} metric on every component of\/
$\Bbb P_\Bbb KV\setminus\S V$.}

\medskip

{\bf Proof.} Let $p\in\Bbb P_\Bbb KV$ be nonisotropic. Let $S$ and $T$
be lifted local fields with $S(p):=s$ and $T(p):=t$.

In order to show that $\big(\nabla_ST-\nabla_TS-[S,T]\big)(p)=0$, we
can assume that the fields $S$ and $T$ are respectively spread from $s$
and $t$ (see Definition 4.1). It follows from Lemma 4.3 that
$\nabla_ST(p)=\nabla_TS(p)=0$. The proof of $[S,T](p)=0$ follows [AGG,
Lemma 4.5.4] : Let $f$ be an smooth function and let $\hat f$ denote
its lift to $V$. By definition,
$\displaystyle T(x)f=\frac d{d\varepsilon}\Big|_{\varepsilon=0}\hat
f\big(x+\pi[x]tx\varepsilon\big)$.
Therefore,
$$S(p)(Tf)=\frac d{d\delta}\Big|_{\delta=0}\Big(\frac
d{d\varepsilon}\Big|_{\varepsilon=0}\hat
f\big(p+sp\delta+\pi[p+sp\delta]t(p+sp\delta)
\varepsilon\big)\Big)=$$
$$=\frac d{d\delta}\Big|_{\delta=0}\Big(\frac
d{d\varepsilon}\Big|_{\varepsilon=0}\hat
f\big(p+sp\delta+\pi[p+sp\delta]tp\varepsilon\big)\Big)=$$
$$=\frac d{d\delta}\Big|_{\delta=0}\bigg(\frac
d{d\varepsilon}\Big|_{\varepsilon=0}\hat
f\Big(p+sp\delta+tp\varepsilon-(p+sp\delta)\frac{k_0\varepsilon\delta}
{1+\delta^2\langle sp,sp\rangle/\langle p,p\rangle}\Big)\bigg),$$
where $k_0:=\langle sp,tp\rangle/\langle p,p\rangle$. Since
$\hat f(pk)=\hat f(p)$ for every $k\in\Bbb K^*$, it follows  that
$$\hat f\big(p(1-k_0\varepsilon\delta)+sp\delta(1-k_0\varepsilon\delta)
+tp\varepsilon\big)=\hat
f\Big(p+sp\delta+tp\frac\varepsilon{1-k_0\varepsilon\delta}\Big).$$
Being $f$ smooth,
$$S(p)(Tf)=\frac d{d\delta}\Big|_{\delta=0}\Big(\frac
d{d\varepsilon}\Big|_{\varepsilon=0}\hat
f\big(p+sp\delta+tp\varepsilon-
(p+sp\delta)k_0\varepsilon\delta\big)\Big)=$$
$$=\frac d{d\delta}\Big|_{\delta=0}\Big(\frac
d{d\varepsilon}\Big|_{\varepsilon=0}\hat
f(p+sp\delta+tp\frac{\varepsilon}{1-k_0\varepsilon\delta})\Big)=\frac
d{d\delta}\Big|_{\delta=0}\Big(\frac
d{d\varepsilon}\Big|_{\varepsilon=0}\hat
f(p+sp\delta+tp\varepsilon)\Big).$$
Hence, $S(p)(Tf)=T(p)(Sf)$, that is, $[S,T](p)=0$.

\smallskip

In order to verify that
$v(S,T)(p)=\big(\nabla_vS(p),T(p)\big)+\big(S(p),\nabla_vT(p)\big)$ for
a tangent vector $v$ at $p$, we~put
$\varphi_1:=\displaystyle\frac
d{d\varepsilon}\Big|_{\varepsilon=0}S(p+vp\varepsilon)$
and
$\varphi_2:=\displaystyle\frac
d{d\varepsilon}\Big|_{\varepsilon=0}T(p+vp\varepsilon)$.
So,
$$\pm\big(\nabla_vS(p),T(p)\big)\dim_\Bbb R\Bbb
K=\pm\big(\pi[p]\varphi_1\pi'[p],T(p)\big)\dim_\Bbb R\Bbb K={\tr}_\Bbb
R\Big(\big(\pi[p]\varphi_1\pi'[p]\big)^*T(p)\Big)={\tr}_\Bbb
R\big(\varphi_1^*T(p)\big),$$
$\pm\big(S(p),\nabla_vT(p)\big)\dim_\Bbb R\Bbb K=\tr_\Bbb R
\big(S^*(p)\varphi_2\big)$, and
$$\pm v(S,T)(p)\dim_\Bbb R\Bbb K=\frac
d{d\varepsilon}\Big|_{\varepsilon=0}{\tr}_\Bbb
R\big(S^*(p+vp\varepsilon)T(p+vp\varepsilon)\big)={\tr}_\Bbb
R\big(\varphi_1^*T(p)\big)+{\tr}_\Bbb R\big(S^*(p)\varphi_2\big).$$
Similar arguments work for the hermitian case
$_\blacksquare$

\medskip

{\bf4.5.~Curvature tensor.} Let $p\in\Bbb P_\Bbb KV$ be nonisotropic
and let $T_1,T_2,S$ be local lifted fields with $T_i(p)=t_i$ and
$S(p)=s$. We wish to express the curvature tensor
$R(T_1,T_2)S(p):=\big(\nabla_{T_2}\nabla_{T_1}S-\nabla_{T_1}
\nabla_{T_2}S+\nabla_{[T_1,T_2]}S\big)(p)$
in terms of the hermitian form. We can assume that the fields $T_i$ and
$S$ are respectively spread from $t_i$ and $s$ (see Definition 4.1). By
Lemma 4.3,
$$\nabla_{T_1}\nabla_{T_2}S(p)=\Big(\frac{d}{d\varepsilon}
\Big|_{\varepsilon=0}\pi[p+t_1p\varepsilon]\big(s\pi[p+t_1p
\varepsilon]t_2-t_2\pi'[p+t_1p\varepsilon]s\big)\pi'[p+t_1p
\varepsilon]\Big)_p.$$
By Lemma 4.2,
$$\Big(\frac{d}{d\varepsilon}\Big|_{\varepsilon=0}\pi[p+t_1p
\varepsilon]s\pi[p+t_1p\varepsilon]t_2\pi'[p+t_1p\varepsilon]
\Big)_p=$$
$$=\big(-(t_1+t_1^*)s\pi[p]t_2\pi'[p]-\pi[p]s(t_1+t_1^*)t_2\pi'[p]+
\pi[p]s\pi[p]t_2(t_1+t_1^*)\big)_p=-st_1^*t_2$$
and
$\displaystyle\Big(\frac{d}{d\varepsilon}\Big|_{\varepsilon=0}
\pi[p+t_1p\varepsilon]t_2\pi'[p+t_1p\varepsilon]s\pi'[p+t_1p
\varepsilon]\Big)_p=t_2t_1^*s$.
In other words, $\nabla_{T_1}\nabla_{T_2}S(p)=-st_1^*t_2-t_2t_1^*s$.
By symmetry, $\nabla_{T_2}\nabla_{T_1}T(p)=-st_2^*t_1-t_1t_2^*s$.
Since $[T_1,T_2](p)=0$ (see the proof of Proposition 4.4), we arrive
at
$$R({t_1},{t_2})s=st_1^*t_2+t_2t_1^*s-st_2^*t_1-t_1t_2^*s\
_\blacksquare$$

{\bf4.6.~Sectional curvature. Constant curvature classic
geometries.} Let $p\in\Bbb P_\Bbb KV$ be nonisotropic. Let
$W\subset\T_p\Bbb P_\Bbb KV$ be a $2$-dimensional $\Bbb R$-vector
subspace such that the metric, being restricted to $W$, is
nondegenerate. The sectional curvature of $W$ is given by
$$SW:=S(t_1,t_2):=\frac{\big(R(t_1,t_2)t_1,t_2\big)}{(t_1,t_1)
(t_2,t_2)-(t_1,t_2)^2}$$
for $\Bbb R$-linearly independent $t_1,t_2\in W$. We can assume that
$t_j=v_j\langle p,-\rangle$ see (2.3), where $v_j\in p^\perp$ and
$\langle v_j,v_j\rangle=\sigma_j\in\{-1,0,+1\}$ for $j=1,2$. In this
way, using the same sign $\pm$ as in (1.4) and applying Remark 2.4, we
obtain
$$\pm(t_1t_1^*t_2,t_2)\dim_\Bbb R\Bbb K={\tr}_\Bbb
R(t_2^*t_1t_1^*t_2)=\dim_\Bbb R\Bbb K\cdot\langle p,p\rangle^2\langle
v_1,v_2\rangle\langle v_2,v_1\rangle.$$
For $k:=\langle v_1,v_2\rangle$, we have
$$\big(R(t_1,t_2)t_1,t_2\big)=\pm\langle
p,p\rangle^2\big(|k|^2+\sigma_1\sigma_2-2\Re(k^2)\big),\quad(t_j,t_j)=
\pm\langle p,p\rangle\sigma_j,\quad(t_1,t_2)=\pm\langle p,p\rangle\Re
k.$$
Hence,
$$SW=\pm\frac{|k|^2+\sigma_1\sigma_2-2\Re(k^2)}{\sigma_1\sigma_2-(\Re
k)^2}=\pm\Big(1+\frac{3|k-\overline k|^2}{4\big(\sigma_1\sigma_2-(\Re
k)^2\big)}\Big),$$
where the last equality follows from the identity
$|k|^2-2\Re(k^2)=\frac34|k-\overline k|^2-(\Re k)^2$. By Lemma 2.1,
$\sigma_1\sigma_2\ne(\Re k)^2$ since $(-,-)$ is nondegenerate over $W$.

\smallskip

Obviously, $SW=\pm1$ if $\Bbb K=\Bbb R$. If $\Bbb K\ne\Bbb R$ and if
$v_1,v_2$ are $\Bbb K$-linearly dependent, then
$\sigma_1\sigma_2=|k|^2$ by Lemma 2.1. In this case,
$|k|=\sigma_1\sigma_2=1$, and it follows from the identity
$|k|^2=|k-\overline k|^2/4+(\Re k)^2$ that $SW=\pm4.$ Since
$v_1,v_2\in p^\perp$ are always $\Bbb K$-linearly dependent if
$\dim_\Bbb KV=2$, we arrive at the

\medskip

{\bf4.7.~Remark.} In every component of $\Bbb P_\Bbb R^n$,
$\Bbb P_\Bbb C^1$, and $\Bbb P_\Bbb H^1$, the sectional curvature is
constant
$_\blacksquare$

\medskip

All the remaining possible values for $SW$ can be extracted from the
above formula. They are displayed in the following table, where
$W=t_1\Bbb R+t_2\Bbb R$, $t_j=v_j\langle p,-\rangle$, and
$v_1,v_2\in p^\perp$ are
$\Bbb K$-linearly independent. The sign $\pm$ is the same as in (1.4).
$$\vbox{\offinterlineskip\hrule\halign{&\vrule#&\strut\quad\hfil#\hfil&
\quad\vrule#&\quad\hfil#\hfil&\quad\vrule#&\quad\hfil#&\quad\vrule#\cr
height3pt&\omit&&\omit&&\omit&\cr&Form over $v_1\Bbb K+v_2\Bbb K$,
$\Bbb K\ne\Bbb R$\hfil&&Metric over $W$&&Sectional curvature\hfill&\cr
height3pt&\omit&&\omit&&\omit&\cr\noalign{\hrule}height3pt&\omit&&\omit
&&\omit&\cr&Indefinite&&Indefinite&&$\pm\,(-\infty,1]$&\cr&Definite&&
Definite&&$\pm\,[1,4)$&\cr&Degenerate&&Definite&&$\pm\,4$&\cr&
Indefinite&&Definite&&$\pm\,(4,\infty)$&\cr
height3pt&\omit&&\omit&&\omit&\cr}\hrule}\ \blacksquare$$

\medskip

\centerline{\bf5.~Parallel transport along geodesics}

\medskip

Let $p\in\Bbb P_\Bbb KV$ be nonisotropic, let $t$ be a tangent vector
at $p$, and let $T$ be the field spread from $t$ (see Definition 4.1).
The smooth (lifted) field
$$\Tn(t)(-):=\frac{T(-)}{\ta(p,-)}$$
is defined out of $\Bbb P_\Bbb Kp^\perp\cup\S V$.

\medskip

{\bf5.1.~Lemma.} {\sl Let\/ $\G$ be a geodesic and let\/ $t$ be a
nonnull tangent vector to\/ $\G$ at a nonisotropic\/ $p\in\G$. Then the
field\/ $\Tn(t)$ is nonnull and tangent to\/ $\G$ wherever defined.}

\medskip

{\bf Proof.} Let $g\in\G$ be nonisotropic and nonorthogonal to $p$.
Clearly, $\varphi:=\Tn(t)(g)\ne0$ since $\pi[g]t\pi'[g]=0$ would imply
$g\in p^\perp$. By Lemma 3.1 (2), $\G=\G W$ with $W=p\Bbb R+tp\Bbb R$.
We can assume that $g\in W$. Hence, $\varphi g\in W$ and $\Tn(t)(g)$ is
tangent to $\G$ at $g$ by Lemma 2.6
$_\blacksquare$

\medskip

{\bf5.2.~Lemma.} {\sl Let\/ $p,q\in\Bbb P_\Bbb KV$ be distinct
nonorthogonal with\/ $p$ nonisotropic. Denote by\/ $\G[p,q]$ the
oriented\/\footnote{In the particular case of a spherical
$\G{\wr}p,q{\wr}$, the segment $\G[p,q]$ is the shortest one from $p$
to $q$.}
segment of the geodesic\/ $\G{\wr}p,q{\wr}$ that does not contain the
point orthogonal to $p$. Let\/ $\varphi:V\to V$ be given by\/
$\varphi=q\langle p,q\rangle^{-1}\langle p,-\rangle$ {\rm(}see\/
{\rm(2.3)}{\rm).} Then\/ $\varphi_p$ is tangent to the\/ {\rm oriented}
segment\/ $\G[p,q]$ at\/ $p$.}

\medskip

{\bf Proof.} The tangent vector $\varphi_p$ does not depend on the
choice of representatives $p,q\in V$. We can assume that
$\langle p,p\rangle=\sigma$ and $\langle p,q\rangle=\sigma a$, where
$\sigma\in\{-1,+1\}$ and $a>0$. Clearly,
$\varphi_p:p\mapsto\pi[p]q(1/a)$. The curve $c_0(t):=p(1-t)+qt$,
$t\in[0,1]$, parameterizes a lift of $\G[p,q]$. Indeed,
$\big\langle p,p(1-t)+qt\big\rangle=0$ means that $(1-a)t=1$, which is
impossible. By Lemma 2.7, the linear map $\dot c(0):p\mapsto\pi[p]q$ is
tangent to $\G[p,q]$ at $p$
$_\blacksquare$

\medskip

{\bf5.3.~Lemma.} {\sl Let\/ $p\in\Bbb P_\Bbb KV$ be nonisotropic, let\/
$t$ be a tangent vector at\/ $p$, and let\/ $T$ be the field spread
from\/ $t$ {\rm(}see Definition\/ {\rm4.1).} Then, for every
nonisotropic\/ $x$,
$$T(x)\big(\ta(p,-)\big)=-2\ta(p,x)\Re\frac{\langle
tx,x\rangle}{\langle x,x\rangle}.$$}
{\bf Proof} is straightforward:
$$T(x)\big(\ta(p,-)\big)=\frac
d{d\varepsilon}\Big|_{\varepsilon=0}\frac{\big\langle
p,x+\pi[x]tx\varepsilon\big\rangle\big\langle
x+\pi[x]tx\varepsilon,p\big\rangle}{\langle p,p\rangle\Big(\langle
x,x\rangle+\varepsilon^2\big\langle\pi[x]tx,\pi[x]tx\big\rangle\Big)}=
\frac{\big\langle p,\pi[x]tx\big\rangle\langle x,p\rangle+\langle
p,x\rangle\big\langle\pi[x]tx,p\big\rangle}{\langle p,p\rangle\langle
x,x\rangle}=$$
$$=-\frac{\langle p,x\rangle\langle x,tx\rangle\langle
x,p\rangle+\langle p,x\rangle\langle tx,x\rangle\langle
x,p\rangle}{\langle p,p\rangle\langle
x,x\rangle^2}=-2\ta(p,x)\Re\frac{\langle tx,x\rangle}{\langle
x,x\rangle}\ _\blacksquare$$

\medskip

{\bf5.4.~Theorem.} {\sl Let\/ $\G$ be a geodesic, let\/ $t$ be a
nonnull tangent vector to\/ $\G$ at a nonisotropic\/ $p\in\G$, and
let\/ $h\in\T_p\L$, where\/ $\L$ stands for the projective line of\/
$\G$. Then, for every nonisotropic\/ $g\in\G$ not orthogonal to\/ $p$,}
$$\nabla_{\Tn(t)(g)}\Tn(h)=0.$$

{\bf Proof.} Denote by $H$ and $T$ the fields respectively spread from
$h$ and $t$ (see Definition 4.1). It suffices to show that
$\Big(\nabla_{T(g)}\displaystyle\frac{H(-)}{\ta(p,-)}\Big)g=0$. By
Lemma 3.1 (2), $\G=\G W$ with $W=p\Bbb R+tp\Bbb R$. We can take
$g\in W$. By Lemmas 4.3 and 5.3,
$$\Big(\nabla_{T(g)}\frac{H(-)}{\ta(p,-)}\Big)g=
T(g)\Big(\frac1{\ta(p,-)}\Big)H(g)g
+\frac1{\ta(p,g)}\big(\nabla_{T(g)}H\big)g=$$
$$=\frac1{\ta(p,g)}\pi[g]\Big(2\frac{\langle tg,g\rangle}{\langle
g,g\rangle}hg+h\pi[g]tg-t\pi'[g]hg\Big).$$
It follows from Lemma 2.6 that $hp=tpk$ for some $k\in\Bbb K$ since
both $h$ and $t$ are tangent to $\L$ at $p$. From
$hp^\perp=tp^\perp=0$, we conclude that $hg=tgk$. Finally, from
$\pi[g]=1-\pi'[g]$, $htg=0$, $\langle tg,g\rangle\in\Bbb R$, and
$hg=tgk$, we obtain
$h\pi[g]tg=-h\pi'[g]tg=-hg\displaystyle\frac{\langle
g,tg\rangle}{\langle g,g\rangle}=-\displaystyle\frac{\langle
tg,g\rangle}{\langle g,g\rangle}hg$
and
$t\pi'[g]hg=t\pi'[g]tgk=tg\displaystyle\frac{\langle
g,tg\rangle}{\langle g,g\rangle}k=\displaystyle\frac{\langle
tg,g\rangle}{\langle g,g\rangle}hg$
$_\blacksquare$

\medskip

Theorem 5.4, Lemma 5.1, and Lemma 3.1 (2) have the following

\medskip

{\bf5.5.~Corollary.} {\sl Out of isotropic points, a geodesic in the
sense of\/ {\rm Example 1.7 (1)} is a geodesic of the Levi-Civita
connection\/ $\nabla$. Every geodesic of this connection appears in
this way}
$_\blacksquare$

\medskip

Of course, Corollary 5.5 can be readily inferred from the standard
characterization of geodesics in symmetric spaces as the trajectory of
certain one-parameter subgroups in the isometry group, but we need
Theorem 5.4 anyway. For example, the theorem provides a formula for the
parallel transport of horizontal vectors along geodesics (see Corollary
5.7).

\medskip

Let $p\in\Bbb P_\Bbb KV$ be nonisotropic, let $t$ be a tangent vector
at $p$, and let $T$ be the field spread from $t$ (see Definition 4.1).
The smooth (lifted) field
$$\Ct(t)(-):=\frac{T(-)}{\sqrt{\ta(p,-)}}$$
is defined at every nonisotropic point in
$\Bbb P_\Bbb KV\setminus\Bbb P_\Bbb Kp^\perp$ that belongs to the
component of $\Bbb P_\Bbb KV\setminus\S V$ containing $p$.

\medskip

{\bf5.6.~Theorem.} {\sl Let\/ $\G$ be a geodesic, let\/ $t$ be a
nonnull tangent vector to\/ $\G$ at a nonisotropic\/ $p\in\G$, and
let\/ $v\in(\T_p\L)^\perp$, where\/ $\L$ stands for the projective line
of\/ $\G$. Then
$$\nabla_{\Tn(t)(g)}\Ct(v)=0$$
for every nonisotropic\/ $g\in\G\setminus\Bbb P_\Bbb Kp^\perp$ that
belongs the component of\/ $\Bbb P_\Bbb KV\setminus\S V$ containing\/
$p$.}

\medskip

{\bf Proof.} Denote by $U$ and $T$ the fields respectively spread from
$v$ and $t$ (see Definition 4.1). It suffices to show that
$\Big(\nabla_{T(g)}\displaystyle\frac{U(-)}{\sqrt{\ta(p,-)}}\Big)g=0$.
By Lemma 3.1 (2), $\G=\G W$ with $W=p\Bbb R+tp\Bbb R$. We can take
$g\in W$. By Lemmas 4.3 and 5.3,
$$\Big(\nabla_{T(g)}\frac{U(-)}{\sqrt{\ta(p,-)}}\Big)g=
T(g)\Big(\frac1{\sqrt{\ta(p,-)}}\Big)U(g)g
+\frac1{\sqrt{\ta(p,g)}}\big(\nabla_{T(g)}U\big)g=$$
$$=\frac1{\sqrt{\ta(p,g)}}\pi[g]\Big(\frac{\langle tg,g\rangle}{\langle
g,g\rangle}vg+v\pi[g]tg-t\pi'[g]vg\Big).$$
By Lemma 2.6, $tpk\langle p,-\rangle\in\T_p\L$ for all $k\in\Bbb K$.
Taking $v\in(\T_p\L)^\perp$ in the form $v=w\langle p,-\rangle$ with
$w\in p^\perp$, we obtain
$\langle p,p\rangle\Re\langle w,tpk\rangle=0$. This implies that
$w\in(p\Bbb K+tp\Bbb K)^\perp$, $vg\in(p\Bbb K+tp\Bbb K)^\perp$,
and $\pi'[g]vg=0$. Finally, as in the proof of Theorem 5.4,
$\displaystyle v\pi[g]tg=-vg\frac{\langle g,tg\rangle}{\langle
g,g\rangle}=-\frac{\langle tg,g\rangle}{\langle g,g\rangle}vg$
$_\blacksquare$

\medskip

Let $\L$ be a noneuclidean projective line and let $p\in\L$ be
nonisotropic. It easily follows from the identification (2.3) that
$\T_p\Bbb P_\Bbb KV=\T_p\L\oplus(\T_p\L)^\perp$. Hence, every tangent
vector $t\in\T_p\Bbb P_\Bbb KV$ decomposes as $t=h+v$, where
$h\in\T_p\L$ and $v\in(\T_p\L)^\perp$. This decomposition is called
{\it horizontal-vertical.} Under the assumption that $\L$ is spanned by
$p$ and $q$, the horizontal-vertical decomposition is
$t=\pi'[w]t+\pi[w]t$, where $w:=\pi[p]q$.

\medskip

{\bf5.7.~Corollary.} {\sl Let\/ $\L$ be a noneuclidean projective line
spanned by distinct, nonisotropic, and nonorthogonal points\/
$p,q\in\Bbb P_\Bbb KV$ of the same signature. Let\/ $t=h+v$ be the
horizontal-vertical decomposition of\/ $t\in\T_p\Bbb P_\Bbb KV$ with
respect to\/ $\L$. Then the parallel transport of\/ $t$ from\/ $p$ to\/
$q$ along\/ $\G[p,q]$ is given by\/ $\Tn\big(h)(q)+\Ct(v)(q)$}
$_\blacksquare$

\medskip

The above corollary expresses the parallel transport along geodesics in
a component of $\Bbb P_\Bbb KV$. However, in particular cases, some
parallel transport can be performed even if the nonisotropic and
nonorthogonal points $p,q$ lie in {\bf different} components of
$\Bbb P_\Bbb KV$ (we just `bat an eye' while passing through $\S V$) :
For a horizontal vector $h$, $\Tn(h)(q)$ gives a parallel transport of
$h$ along $\G[p,q]$. When $\Bbb K=\Bbb C$, for a vertical vector $v$,
$\Ct(v)(q)$ gives a parallel transport of $v$ along $\G[p,q]$ (we fix
the sign of $\sqrt{\ta(p,q)}\in\Bbb Ri$).

\smallskip

It remains to study the parallel transport along euclidean geodesics.
Let $p\in\Bbb P_\Bbb KV$ be nonisotropic, let $s$ be a tangent vector
at $p$, and let $S$ be the field spread from $s$ (see Definition 4.1).
The smooth (lifted) vector field
$$\Eu(s)(x):=\frac12\big(\pi[p]\pi'[x]s\big)_x+S(x)$$
is defined out of isotropic points. Clearly, $\Eu(s)(p)=S(p)=s$.

\medskip

{\bf5.8.~Theorem.} {\sl Let\/ $\G$ be an euclidean geodesic, let\/ $t$
be a nonnull tangent vector to\/ $\G$ at a nonisotropic\/ $p\in\G$, and
let\/ $s\in\T_p\Bbb P_\Bbb KV$. Then, for every nonisotropic\/
$g\in\G$,}
$$\nabla_{\Tn(t)(g)}\Eu(s)=0.$$

{\bf Proof.} It suffices to show that
$\big(\nabla_{T(g)}\Eu(s)\big)g=0$, where $T$ is the field spread from
$t$ (see Definition~4.1). By Lemma 3.1 (2), $\G=\G W$ with
$W=p\Bbb R+tp\Bbb R$. We can take $g\in W$. Note that, being orthogonal
to $p$, each one of $tp$, $tg$, and $\pi[p]g$ represents the only
isotropic point $u\in\G$. Clearly, $\langle u,\G\rangle=0$. It follows
that $\pi[g]t=\pi[p]t=t$. Hence, $s\pi[g]t=st=0$. Also,
$\pi'[g]\pi[p]\pi'[g]=0$. Now, using $\pi[g](t_g)^*=(t_g)^*g=0$, we
obtain
$$2\big(\nabla_{T(g)}\Eu(s)\big)g=\pi[g]\Big(\frac
d{d\varepsilon}\Big|_{\varepsilon=0}\pi[g+t_gg\varepsilon]
\pi[p]\pi'[g+t_gg\varepsilon]s\pi'[g+t_gg\varepsilon]\Big)g+2\pi[g]
s\pi[g]tg-2\pi[g]t\pi'[g]sg=$$
$$=-\pi[g]\big(t_g+(t_g)^*\big)\pi[p]\pi'[g]sg+\pi[g]
\pi[p]\big(t_g+(t_g)^*\big)sg+\pi[g]\pi[p]\pi'[g]s
\big(t_g+(t_g)^*\big)g-2\pi[g]t\pi'[g]sg=$$
$$=\pi[g]\pi[p]\big(t_g+(t_g)^*\big)sg-2\pi[g]t\pi'[g]sg$$
by Lemmas 4.2 and 4.3. Since $(\varphi\psi)^*=\psi^*\varphi^*$,
$\langle g,t^*g\rangle=\langle tg,g\rangle=0$, $\pi[g]\pi[p]g=\pi[p]g$,
and the projections are self-adjoint, we obtain
$$\pi[g]\pi[p](t_g)^*sg=\pi[g]\pi[p]\pi'[g]t^*\pi[g]sg=
\pi[g]\pi[p]\pi'[g]\Big(t^*sg-t^*g\frac{\langle g,sg\rangle}{\langle
g,g\rangle}\Big)=$$
$$=\pi[g]\pi[p]\Big(g\frac{\langle g,t^*sg\rangle}{\langle
g,g\rangle}-g\frac{\langle g,t^*g\rangle\langle g,sg\rangle}{\langle
g,g\rangle^2}\Big)=\pi[p]g\frac{\langle tg,sg\rangle}{\langle
g,g\rangle}.$$
It follows from $\pi[p]t=\pi[g]t=t$ and $sg\in p^\perp$ that
$$\pi[g]\pi[p]t_gsg=\pi[g]t\pi'[g]sg=tg\displaystyle\frac{\langle
g,sg\rangle}{\langle
g,g\rangle}=tg\frac{\big\langle\pi[p]g,sg\big\rangle}{\langle
g,g\rangle}.$$
It remains to observe that $\pi[p]g$ and $tg$ are $\Bbb R$-proportional
$_\blacksquare$

\medskip

{\bf5.9.~Corollary.} {\sl Let\/ $p,q\in\Bbb P_\Bbb KV$ be distinct and
nonisotropic points that span an euclidean projective line and let\/
$t\in\T_p\Bbb P_\Bbb KV$. Then the parallel transport of\/ $t$ from\/
$p$ to\/ $q$ along\/ $\G[p,q]$ is given by\/
$\Eu(t)(q)$}~
$_\blacksquare$

\bigskip

\centerline{\bf6.~Complex hyperbolic examples}

\medskip

The three examples below concern complex hyperbolic geometry. For basic
background on the subject, see [Gol] or [AGG, Section 4]. As in Example
1.6 (4), we take $\Bbb K=\Bbb C$, $\dim_\Bbb CV=3$, the form of
signature $++-$ and the sign $-$ in the definition (1.5) of the
hermitian metric. Thus, $\B V$ is the complex hyperbolic plane
$\Bbb H^2_\Bbb C$.

\medskip

{\bf6.1.~Example: area formula.} Let $p_1,p_2,p_3\in\B V\cup\S V$ be
points in a complex geodesic $\L$. With the use of vertical parallel
transport, we will show that the oriented area of the plane triangle
$\triangle(p_1,p_2,p_3)$ is given by
$$\Area\triangle(p_1,p_2,p_3)=-\textstyle\frac{1}{2}\arg\big(-\langle
p_1,p_2\rangle\langle p_2,p_3\rangle\langle
p_3,p_1\rangle\big),\leqno{\bold{(6.2)}}$$
where $\arg$ varies in $[-\pi,\pi]$.

First, we take $p_j\notin\S V$, $j=1,2,3$. We have
$\L=\Bbb P_\Bbb Cp^\perp$, where $p\in\E V$ is the polar point to $\L$
(for the definition of polar point, see the beginning of Example 3.6 or
[AGG, Subsection 4.1.6]). By Lemma 2.6,
$(\T_q\L)^\perp=p\Bbb C\langle q,-\rangle$ for every
$q\in\L\setminus\S V$.
Let $v:=pc\langle p_1,-\rangle\in(\T_{p_1}\L)^\perp$, $c\in\Bbb C^*$.
Making the parallel transport of $v$ along the segment of geodesic
$\G[p_1,p_2]$, then along $\G[p_2,p_3]$, and finally along
$\G[p_3,p_1]$, we end up with some $v'\in(\T_{p_1}\L)^\perp$.
By Corollary 5.7,
$$v'=\frac{\pi[p_1]\pi[p_3]\pi[p_2]v\pi'[p_2]\pi'[p_3]\pi'[p_1]}
{\sqrt{\ta(p_1,p_2)\ta(p_2,p_3)\ta(p_3,p_1)}}=\frac{pc\langle
p_1,p_2\rangle\langle p_2,p_3\rangle\langle p_3,p_1\rangle\langle
p_1,-\rangle}{\langle p_2,p_2\rangle\langle p_3,p_3\rangle\langle
p_1,p_1\rangle\sqrt{\ta(p_1,p_2)\ta(p_2,p_3)\ta(p_3,p_1)}}$$
because $p\in p_j^\perp$. Clearly, $(\T_{p_1}\L)^\perp$ is a
one-dimensional $\Bbb C$-vector space. The oriented angle
$\angle(v,v')$ from $v$ to $v'$, taken in $[-\pi,\pi]$, is an additive
measure of a triangle. Hence, it is proportional to the oriented area
of $\triangle(p_1,p_2,p_3)$. In terms of the hermitian metric (1.5),
$$\angle(v,v')=\arg\langle v,v'\rangle=\arg\big(-\langle
p_1,p_2\rangle\langle p_2,p_3\rangle\langle p_3,p_1\rangle\big)$$
due to $p\in\E V$ and $p_2,p_3\in\B V$. The formula is extendable to
isotropic points. Considering a suitable ideal triangle, we find the
factor of proportionality $-1/2$ in (6.2).

\smallskip

The obtained formula (without orientation taken into account) can be
found in [Gol]. Using the horizontal parallel transport instead of
the vertical one, we would arrive at the well-known area formula in
terms of the angles. Curiously, the formula (6.2) seems to appear more
naturally in the context of complex hyperbolic geometry. A similar
formula holds for a plane spherical triangle
$_\blacksquare$

\medskip

{\bf6.3.~Example: some geometry behind the angle between
bisectors.} Let $B_1$ and $B_2$ be bisectors in $\Bbb H^2_\Bbb C$ with
hyperbolic real spines $\G_1$ and $\G_2$. Assume that these bisectors
share a common slice $S$ whose polar point is $p\in\E V$. Let
$v_j\in\S V\cap\G_j$ denote some vertex of $B_j$, $j=1,2$. Then the
point $q_j:=\pi[p]v_j$ is the intersection point of the real spine of
$B_j$ with the slice $S$. Denote by $\G[q_j,v_j)\subset\G_j$ the
oriented segment of the real spine that starts with $q_j$ and ends with
$v_j$. Let $B[q_j,v_j)\subset B_j$ denote the corresponding oriented
segment of bisector: $B[q_j,v_j)$ is oriented with respect to the
orientation of $\G[q_j,v_j)$ and to the natural orientation of its
slices. Define
$$u:=1-\frac{\langle v_2,v_1\rangle\langle p,p\rangle}{\langle
v_2,p\rangle\langle p,v_1\rangle}.$$
In other words, $u=1-\displaystyle\frac1{\eta(v_1,v_2,p)}$, where
$\eta(v_1,v_2,p)$ is {\it Goldman's invariant\/} [Gol].

Let $q\in S$. We choose representatives $p,v_1,v_2\in V$ such that
$\langle p,p\rangle=\langle p,v_j\rangle=1$. Thus,
$$q_j=v_j-p,\quad\langle q_j,v_j\rangle=-1,\quad\langle
q_j,q\rangle=\langle v_j,q\rangle,\quad\langle q_j,q_j\rangle=-1,$$
$$\pi[q_j]v_j=p,\quad\langle q_2,q_1\rangle=\langle
v_2,v_1\rangle-1=-u,\quad\ta(q_1,q_2)=|u|^2.$$
In particular, $u\ne0$. According to [AGG, Proposition 4.2.11 and Lemma
4.2.15],
$$n(q,q_j,v_j)=\Big(q_j\frac{\langle v_j,q\rangle}{\langle
v_j,q_j\rangle}-v_j\frac{\langle q_j,q\rangle}{\langle
q_j,v_j\rangle}\Big)i\,\langle q,-\rangle=p\langle
v_j,q\rangle i\,\langle q,-\rangle$$
is a normal vector to the oriented segment $B[q_j,v_j)$ at $q$. Both
normal vectors in question belong to the $\Bbb C$-vector space
$(\T_qS)^\perp$ and, therefore, the oriented angle
$\angle\big(q,B[q_1,v_1),B[q_2,v_2)\big)$ from $B[q_1,v_1)$ to
$B[q_2,v_2)$ at $q$ can be calculated as
$$\angle\big(q,B[q_1,v_1),B[q_2,v_2)\big)=\arg\big\langle
n(q,q_1,v_1),n(q,q_2,v_2)\big\rangle=\arg\big(-\langle
q,q\rangle\langle q,v_1\rangle\langle v_2,q\rangle\big)=$$
$$=\arg\big(\langle q,v_1\rangle\langle
v_2,q\rangle\big)=\arg\big(\langle q,q_1\rangle\langle
q_2,q\rangle\big)=\arg\big(-u\langle q,q_1\rangle\langle
q_1,q_2\rangle\langle q_2,q\rangle\big)$$
since $-u\langle q_1,q_2\rangle=|u|^2$. In other words, using the
previous example,
$$\angle\big(q,B[q_1,v_1),B[q_2,v_2)\big)\equiv\arg
u-2\Area\Delta(q,q_1,q_2)\mod2\pi.$$

We can see that the angle in question is composed of two parts. The
{\it constant angle\/} $\arg u$ is independent of $q\in S$ (in [Hsi3],
this angle is called prespinal). The {\it nonconstant angle\/}
$-2\Area(q,q_1,q_2)$ depends only on the mutual position of $q,q_1,q_2$
in $S$. Let us show that the constant angle is the angle from the real
spine $\G[q_1,v_1)$ to the real spine $\G[q_2,v_2)$ measured with the
help of parallel transport along the segment of geodesic $\G[q_1,q_2]$.

By Lemma 5.2,
$t_j:=\pi[q_j]v_j\langle q_j,v_j\rangle^{-1}\langle
q_j,-\rangle=-p\langle q_j,-\rangle$
is tangent to $\G[q_j,v_j)$ at $q_j$. By Corollary~5.7, the parallel
transport of $t_1$ along $\G[q_1,q_2]$ is given by
$$\Ct(t_1)(q_2)=\frac{\pi[q_2]t_1\pi'[q_2]}{\sqrt{\ta(q_1,q_2)}}=
-\frac{\pi[q_2]p\langle q_1,q_2\rangle\langle
q_2,-\rangle}{|u|\langle
q_2,q_2\rangle}=-\frac{\overline u}{|u|}p\langle
q_2,-\rangle=\frac{\overline u}{|u|}t_2.$$
This implies the result, illustrated by the following picture:

\medskip

\noindent
$\hskip72pt\vcenter
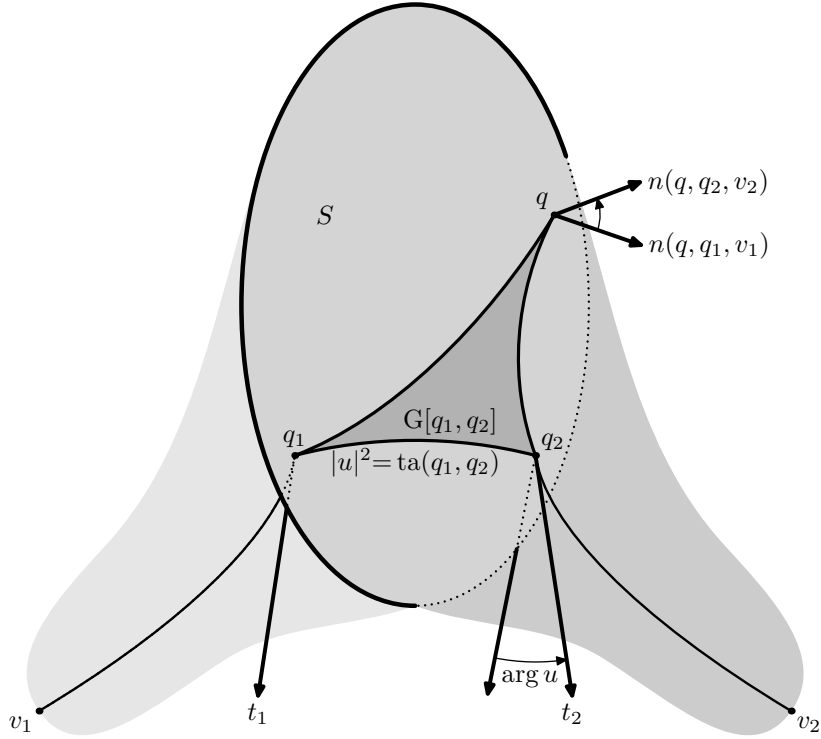
{\hbox{\epsfbox{picture1.eps}}}$

\vskip5pt

It easily follows from Sylvester's criterion that $u$ completely
characterizes the configuration of $B[q_1,v_1)$ and $B[q_2,v_2)$ and
that every $u\in\Bbb C$ with $|u|\ge1$ is possible. The geometric
meaning of $u$ is clear now: $|u|^2$ is the tance between the complex
spines of the bisectors and $\arg u$ is the angle between their real
spines, in the above sense
$_\blacksquare$

\medskip

{\bf6.4.~Example: meridional and parallel transports.} Let $B$ be a
bisector in $\Bbb P_\Bbb CV$ as introduced in Example 1.7 (4), let $\G$
and $\L$ be the real and complex spines of $B$, and let $p_1,p_2\in\G$
be distinct, nonisotropic, and nonorthogonal points. Denote by $S_j$
the slice of $B$ that contains $p_j$, $j=1,2$. Take $q_1\in S_1$
different from the focus $f$ of $B$. The slice $S_j$ is spanned by
$p_j$ and $f$. By Lemma 2.6, the complex spine and the slices are
orthogonal.

The vector
$v:=\pi[p_1]q_1\langle p_1,q_1\rangle^{-1}\langle p_1,-\rangle$ is
tangent to $\G[p_1,q_1]\subset\S_1$ at $p_1$ by Lemma 5.2 and is thus
orthogonal to the complex spine of $B$. Let $\Ct(v)(p_2)$ denote the
parallel transport of $v$ from $p_1$ to $p_2$ along $\G[p_1,p_2]$ given
by Corollary 5.7 and by the considerations right after it. Then there
exists a unique $q_2\in S_2$ such that
$$\pi[p_2]q_2\langle p_2,q_2\rangle^{-1}\langle
p_2,-\rangle=\Ct(v)(p_2).$$
(This can be seen by considering $q_2$ in the form $q_2=p_2+fc$,
$c\in\Bbb C$.) We call $q_2$ the {\it meridional transport\/} of $q_1$
from $p_1$ to $p_2$ along $\G[p_1,p_2]$. In explicit terms,
$$q_2=p_2\langle p_1,q_1\rangle\sqrt{\ta(p_1,p_2)}+\pi[p_1]q_1\langle
p_1,p_2\rangle.$$
The meridional transport identifies almost all slices of the bisector
(the only exceptions are the slices tangent to $\S V$, if they exist).
Such identification, called the {\it slice identification,} is an
important tool for constructing and characterizing complex hyperbolic
manifolds in [AGG].

The meridional and parallel transports are related as follows. As is
easy to see, every slice $S$ of $B$ has the form
$S=\Bbb P_\Bbb Cg^\perp$, where $g\in\G$ is the polar point to $S$. If
$g$ is nonisotropic, we associate to every nonnull tangent vector
$t\in\T_g\Bbb P_\Bbb CV$ the point $tg\in S$. Denote by $g_j\in\G$ the
polar points to $S_j$. The parallel transport along $\G[g_1,g_2]$
produces the meridional transport of the associated points:

\vskip3pt

\noindent
\centerline{$\vcenter
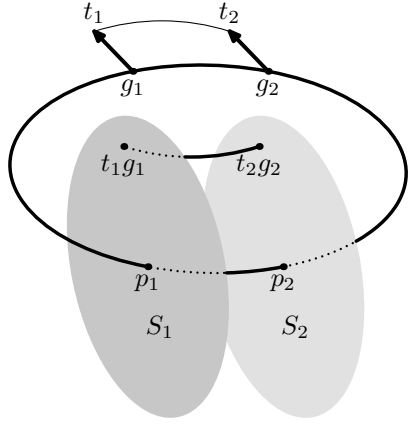
{\hbox{\epsfbox{picture2.eps}}}$}

\vskip3pt

\noindent
Indeed, $g_1,g_2$ are nonorthogonal and nonisotropic. Let $t_1$ be a
tangent vector at $g_1$. By Corollary 5.7, the parallel transport of
$t_1$ from $g_1$ to $g_2$ along $\G[g_1,g_2]$ is given by
$$t_2:=\Tn(h)(g_2)+\Ct(v)(g_2)=\Big(\frac h{\ta(g_1,g_2)}+\frac
v{\sqrt{\ta(g_1,g_2)}}\Big)_{g_2},$$
where $t_1=h+v$ is the horizontal-vertical decomposition of $t_1$ with
respect to $\L$, that is, $h\in\T_{g_1}\L$ and
$v\in(\T_{g_1}\L)^\perp$. We can assume that $h\ne0$ (otherwise, the
focus $f$ is the point associated to both $t_1$ and $t_2$). It is easy
to see that $\ta(g_1,g_2)=\ta(p_1,p_2)$. Since $\pi'[g_1]g_2$ and $g_1$
are $\Bbb C^*$-proportional, the point in $S_2$ associated to $t_2$ has
the form
$$t_2g_2=\frac{\pi[g_2]hg_2}{\ta(g_1,g_2)}+
\frac{\pi[g_2]vg_2}{\sqrt{\ta(g_1,g_2)}}\sim\frac{\pi[g_2]hg_1\langle
p_1,p_1\rangle\langle p_2,p_2\rangle}{\langle
p_2,p_1\rangle}\sqrt{\ta(p_1,p_2)}+vg_1\langle p_1,p_2\rangle,$$
where $\sim$ means $\Bbb C^*$-proportionality. By Lemma 2.6,
$hg_1\in(p_1\Bbb C+g_1\Bbb C)\cap g_1^\perp=p_1\Bbb C$ because
$h\in\T_{g_1}\L$. Also, $vg_1\in f\Bbb C$. From $t_1=h+v$ and from the
orthogonal decomposition $p_2\Bbb C+g_2\Bbb C$, it follows now that
$\pi[p_1]t_1g_1=vg_1$ and
$\pi[g_2]hg_1=\pi'[p_2]hg_1=p_2\displaystyle\frac{\langle
p_2,hg_1\rangle}{\langle p_2,p_2\rangle}$.
It remains to observe that $hg_1\in p_1\Bbb C$ implies that
$\langle p_2,hg_1\rangle=\big\langle\pi'[p_1]p_2,hg_1\big\rangle=
\displaystyle\frac{\langle p_1,hg_1\rangle\langle
p_2,p_1\rangle}{\langle p_1,p_1\rangle}$
$_\blacksquare$

\bigskip

{\bf Acknowledgement.} We are very grateful to the referees for
valuable remarks.

\newpage

\centerline{\bf References}

\medskip

[AGG] S.~Anan$'$in, C.~H.~Grossi, N.~Gusevskii, {\it Complex hyperbolic
structures on disc bundles over surfaces,} to appear in
Int.~Math.~Res.~Not., see also http://arxiv.org/abs/math/0511741

[AGK] D.~V.~Alekseevsky, B.~Guilfoyle, W.~Klingenberg, {\it On the
geometry of spaces of oriented geodesics,} to appear in
Ann.~Glob.~Anal.~Geom., see also http://arxiv.org/abs/0911.2602

[AGoG] S.~Anan$'$in, E.~C.~B.~Gon\c calves, C.~H.~Grossi, {\it
Grassmannians and conformal structure on absolutes,} preprint
http://arxiv.org/abs/0907.4469

[AGr] S.~Anan$'$in, C.~H.~Grossi, {\it Differential geometry of
grassmannians and Pl\"ucker map,} preprint
http://arxiv.org/abs/0907.4470

[Arn1] V.~I.~Arnold, {\it Mathematical Methods of Classical Mechanics,}
GTM {\bf60}, Springer-Verlag, New York, 1997. xx+509 pp.

[Arn2] V.~I.~Arnold, {\it Lobachevsky triangle altitudes theorem as the
Jacobi identity in the Lie algebra of quadratic forms on symplectic
plane,} J.~Geom.~Phys.~{\bf53} (2005), no.~4, 421--427

[BeP] R.~Benedetti, C.~Petronio, {\it Lectures on Hyperbolic Geometry,}
Universitext, Springer-Verlag, 2008. ix+323 pp.

[GeG] N.~Georgiou, B.~Guilfoyle, {\it On the space of oriented
geodesics of hyperbolic\/ $3$-space,} Rocky Mountain J.~Math.~{\bf40}
(2010), no.~4, 1183--1219

[ChG] S.~S.~Chen, L.~Greenberg, {\it Hyperbolic Spaces,} Contributions
to analysis (a collection of papers dedicated to Lipman Bers), Academic
Press (1974), 49--87

[ChK] Y.~Cho, H.~Kim, {\it The analytic continuation of hyperbolic
space,} preprint\newline http://arxiv.org/abs/math/0612372

[Gir] G.~Giraud, {\it Sur certaines fonctions automorphes de deux
variables,} Ann.~Ec.~Norm.~(3), {\bf38} (1921), 43--164

[Gol] W.~M.~Goldman, {\it Complex Hyperbolic Geometry,} Oxford
Mathematical Monographs. Oxford Science Publications. The Clarendon
Press, Oxford University Press, New York, 1999. xx+316 pp.

[Gro] C.~H.~Grossi, {\it Elementary tools for classic and complex
hyperbolic geometries,} PhD thesis, State University of Campinas,
20 September 2006, iii+134 pp.

[GuK1] B.~Guilfoyle, W.~Klingenberg, {\it An indefinite K\"ahler metric
on the space of oriented lines,} J.~London Math.~Soc.~{\bf72} (2005),
497–-509

[GuK2] B.~Guilfoyle, W.~Klingenberg, {\it Proof of the Caratheodory
conjecture by mean curvature flow in the space of oriented affine
lines,} preprint http://arxiv.org/abs/0808.0851

[How] H.~Roger, {\it Remarks on classical invariant theory,}
Trans.~Amer.~Math.~Soc.~{\bf313} (1989), no.~2, 539--570

[HSa] J.~Hakim, H.~Sandler, {\it Application of Bruhat decomposition to
complex hyperbolic geometry,} J.~Geom.~Anal.~{\bf10} (2000), 435--453

[Hsi1] P.~H.~Hsieh, {\it Linear submanifolds and bisectors in\/
$\Bbb CH^n$,} Forum.~Math.~{\bf10} (1998), 413--434

[Hsi2] P.~H.~Hsieh, {\it Semilinear submanifolds in complex hyperbolic
space,} Forum.~Math.~{\bf11} (1999), 673--694

[Hsi3] P.~H.~Hsieh, {\it Cotranchal bisectors in complex hyperbolic
space,} Geometriae Dedicata {\bf97} (2003), 93--98

[Kle] F.~Klein, {\it Vorlesungen \"uber h\"ohere Geometrie,}
Grundlehren der mathematischen Wissenschaften {\bf22}, Julius Springer,
Berlin, 1926. viii+405 pp.

[Lan] S.~Lang, {\it Algebra,} Advanced Book Program. Addison-Wesley
Publishing Company, Inc., California, 1984. xx+714 pp.

[Man] Yu.~I.~Manin, {\it Gauge Field Theory and Complex Geometry,} A
Series of Comprehensive Studies in Mathematics {\bf289},
Springer-Verlag, 1988. x+297 pp.

[Mos] G.~D.~Mostow, {\it On a remarkable class of polyhedra in complex
hyperbolic space,} Pacific J.~Math. {\bf86} (1980), 171--276

[Sal1] M.~Salvai, {\it On the geometry of the space of oriented lines
in Euclidean space,} Manuscripta Math.~{\bf118} (2005), 181--189

[Sal2] M.~Salvai, {\it Geometry of the space of oriented lines in
hyperbolic space,} Glasgow Math.~J.~{\bf49} (2007), 357--366

[San] H.~Sandler, {\it Distance formulas in complex hyperbolic space,}
Forum Math.~{\bf8} (1996), no.~1, 93--106

[Stu] E.~Study, {\it Von den Bewegungen und Umlegungen\/ {\rm I,}
{\rm II,}} Math.~Ann.~{\bf34} (1891), 441--566

[Thu] W.~P.~Thurston, {\it Three-Dimensional Geometry and
Topology.~{\rm I,}} Princeton Mathematical Series, Princeton University
Press, Princeton, 1997. x+311 pp.

[Wey] H.~Weyl, {\it The classical groups,} Princeton Mathematical
Series, Princeton University Press, Princeton, 1939. xiv+320 pp.

[Wil1] N.~J.~Wildberger, {\it Universal Hyperbolic Geometry\/ {\rm I:}
Trigonometry,} preprint\newline http://arxiv.org/abs/0909.1377

[Wil2] N.~J.~Wildberger, {\it Universal Hyperbolic Geometry\/ {\rm
II:} A pictorial overview,} KoG {\bf 14} (2010), 1--23, see also
http://arxiv.org/abs/1012.0880

\enddocument